\documentclass{birkjour}
\usepackage{amsmath}
\usepackage{amsfonts}
\usepackage{amssymb}
\usepackage{float}
\usepackage{graphicx}
\usepackage{epstopdf}
\usepackage{hyperref}
\usepackage{multirow}
\usepackage[misc]{ifsym}
\providecommand{\U}[1]{\protect\rule{.1in}{.1in}}
\providecommand{\U}[1]{\protect\rule{.1in}{.1in}}
\providecommand{\U}[1]{\protect\rule{.1in}{.1in}}
\providecommand{\U}[1]{\protect\rule{.1in}{.1in}}

\newtheorem{theorem}{Theorem}

\newtheorem{algorithm}[theorem]{Algorithm}

\newtheorem{condition}[theorem]{Condition}

\numberwithin{equation}{section} \numberwithin{theorem}{section}
\newtheorem{definition}[theorem]{Definition}
\newtheorem{example}[theorem]{Example}

\newtheorem{lemma}[theorem]{Lemma}

\newtheorem{remark}[theorem]{Remark}

\begin{document}

\title[Projection and contraction algorithms with outer perturbations]{Convergence of projection and contraction algorithms with outer perturbations and their applications to sparse signals recovery}

\author[Q.-L. Dong]{Qiao-Li Dong}
\address{Tianjin Key Laboratory for Advanced Signal Processing, College
of Science, Civil Aviation University of China, Tianjin 300300, China}

\author[A. Gibali]{Aviv Gibali}
\address{Department of Mathematics, ORT Braude College, 2161002 Karmiel,
Israel}
\email{\Letter\hspace{.01in} avivg@braude.ac.il}

\author[D. Jiang]{Dan Jiang}
\address{Tianjin Key Laboratory for Advanced Signal Processing, College
of Science, Civil Aviation University of China, Tianjin 300300, China}

\author[S.-H. Ke]{Shang-Hong Ke}
\address{Tianjin Key Laboratory for Advanced Signal Processing, College
of Science, Civil Aviation University of China, Tianjin 300300, China}



\date{Submitted: April 24, 2017. Revised: August 20, 2017. Accepted for publication in \textit{Journal of Fixed Point Theory and Applications}}

\begin{abstract}
In this paper we study the bounded perturbation
resilience of projection and contraction algorithms
for solving variational inequality (VI) problems in real Hilbert spaces.
Under typical and standard assumptions of monotonicity and Lipschitz continuity
of the VI's associated mapping, convergence of the perturbed projection and contraction algorithms is proved.
Based on the bounded perturbed resilience of projection and contraction algorithms,
we present some inertial projection and contraction algorithms.
In addition we show that the perturbed algorithms converges at the rate of $O(1/t)$.
\end{abstract}

\maketitle

\section{Introduction}

In this article, we are concerned with the classical variational inequality (VI) problem, which is to find
a point $x^*\in C$ such that
\begin{equation}
\label{vip}
\langle F(x^*),x-x^* \rangle \geq 0 ,\,\,\,\,\forall x \in C.
\end{equation}
Where $C$ is a closed convex set in Hilbert space $\mathcal{H}$, $\langle \cdot,\cdot \rangle$ denotes the inner product in $\mathcal{H}$ and $F: \mathcal{H} \rightarrow \mathcal{H}$ is the VI associated mapping.

This problem is a fundamental problem in optimization theory and related fields. It captures various applications, such as
partial differential equations, optimal control, and mathematical programming.
There exist many iterative algorithms for solving the VI (\ref{vip}); For example the \textit{extragradient method} of Korpelevich \cite{Korpelevich} (also Antipin \cite{Antipin}), in which at each iteration of the algorithm, in order to get the next iterate $x^{k+1}$, two orthogonal projections onto $C$ are calculated, according to the following iterative step. Given the current
iterate $x^{k},$ calculate%
\begin{equation}
\left\{
\begin{array}{l}
y^{k}=P_{C}(x^{k}-\beta_k F(x^{k}))\medskip \\
x^{k+1}=P_{C}(x^{k}-\beta_k F(y^{k}))%
\end{array}%
\right.  \label{EG-1}
\end{equation}%
where $\beta_k\in (0,1/L)$, and $L$ is the Lipschitz constant of $F$ (or $\beta_k$ is replaced by a sequence of $\{\beta_k\}^\infty_{k=1}$ which is updated by some adaptive procedure). For an extensive and excellent book on theory, algorithms and applications of VIs see Facchinei and Pang book, \cite{PF03}. In this matter see also the comparative numerical study regarding gradient and extragradient methods for solving VIs \cite{gjkt16}.

In this paper we wish to focus on a close but different type of algorithms, known as projection and contraction algorithms (PC-algorithms). They are called projection and contraction algorithms, according to \cite{He}, because in each iteration projections are used and the distance of the iterates to the solution set of the VI monotonically converges to zero.

He \cite{He} and Sun \cite{Sun} developed a projection and contraction algorithm, which consist of
two steps.
The first one produces the $k$-th iterate point $y^k$ in the the same way as in the extragradent method:
\begin{equation}
\label{PC-y}
y^k=P_C(x^k-\beta_k F(x^k)).
\end{equation}
but the second update of the next iteration $x^{k+1}$ step is updated via the following PC-algorithms:
\begin{equation}
\label{PC-x1}
(\hbox{PC-algorithm I})\quad x^{k+1}=x^k-\gamma \varrho_k d(x^k,y^k)
\end{equation}
or
\begin{equation}
\label{PC-x2}
(\hbox{PC-algorithm II})\quad x^{k+1}=P_C(x^k-\gamma\varrho_k\beta_k F(y^k)),
\end{equation}
where $\gamma\in (0,2),$  $\beta_k\in (0,1/L)$ (or $\{\beta_k\}$ which is updated by some self-adaptive rule),
\begin{equation}
\label{PC-d}
d(x^k,y^k):=(x^k-y^k)-\beta_k (F(x^k)-F(y^k))
\end{equation}
and
\begin{equation}
 \label{pc-varrho}
\varrho_k:=\frac{\langle x^k-y^k,d(x^k,y^k)\rangle}{\|d(x^k,y^k)\|^2}.
\end{equation}

Cai {\it et al.} \cite[Theorem 4.1]{He1} proved the convergence of the PC-algorithms in Euclidean spaces. Dong {\it et al.} \cite[Theorem 3.1]{Dong1} extended the results of \cite{He1} to Hilbert spaces and proved the weak convergence of the PC-algorithm (\ref{PC-x2}).
In order to present a direct consequence from these two results we need to assume the following conditions on the VI (\ref{vip}).

\begin{condition}
\label{con:Condition 1.1}
The solution set of (\ref{vip}), denoted by $SOL(C,F)$, is nonempty.
\end{condition}

\begin{condition}
\label{con:Condition 1.2}
The mapping $F$ is \texttt{monotone}, i.e.,
\begin{equation}
\label{monotone}
\langle F(x)-F(y),x-y\rangle\geq0\,\,\,\,\forall x,y\in \mathcal{H}.
\end{equation}
\end{condition}

\begin{condition}
\label{con:Condition 1.3}
The mapping $F$ is \texttt{Lipschitz-continuous} on $\mathcal{H}$ with constant $L>0$, i.e., there exists $L>0$ such that
\begin{equation}
\label{Lipschitz}
\|F(x)-F(y)\|\leq L\|x-y\|\,\,\,\,\forall x,y\in \mathcal{H}.
\end{equation}
\end{condition}

Hence, we can now establish the following Theorem derived from \cite{He1} and \cite{Dong1}.

\begin{theorem}
Assume that Conditions \ref{con:Condition 1.1}--\ref{con:Condition 1.3} hold. Then any sequence
$\{x^k\}_{k=0}^\infty$  generated by the projection and contraction algorithms (\ref{PC-y})--(\ref{pc-varrho}) weakly converges
to a solution of the  variational inequality (\ref{vip}).
\end{theorem}

The purpose of this paper is then to prove the bounded perturbation
resilience of the PC-algorithms for solving variational inequality (VI) problem in real Hilbert spaces.
This would enable to apply the Superiorization methodology and also introduce inertial PC-algorithms. Moreover, we show that the perturbed algorithms converge at the rate of $O(1/t)$.

The outline of the paper is a s follows. In Section \ref{sec:Preliminaries} we
present definitions and notions that will be need for the rest of the paper. In
Section \ref{sec:Convergence} the PC-algorithms with outer perturbations are presented and analyzed. Later in Section
\ref{sec:BPR} the bounded perturbation resilience of the PC-algorithms is proved, then in Section \ref{sec:iner-PC}
we construct the inertial PC-algorithms. Finally in Section \ref{sec:num} we compare and demonstrate the algorithms performances with respect to the problem of sparse signal recovery .

\section{Preliminaries}\label{sec:Preliminaries}
Let $\mathcal{H}$ be a real Hilbert space with inner product $\langle \cdot
,\cdot \rangle $ and the induced norm $\Vert \cdot \Vert $, and let $D$ be a
nonempty, closed and convex subset of $\mathcal{H}$. We write $%
x^{k}\rightharpoonup x$ to indicate that the sequence $\left\{ x^{k}\right\}
_{k=0}^{\infty }$ converges weakly to $x$ and $x^{k}\rightarrow x$ to
indicate that the sequence $\left\{ x^{k}\right\} _{k=0}^{\infty }$
converges strongly to $x.$ Given a sequence $\left\{ x^{k}\right\}
_{k=0}^{\infty }$, denote by $\omega _{w}(x^{k})$ its weak $\omega $-limit
set, that is, any $x\in \omega _{w}(x^{k})$ such that there exsists a subsequence $%
\left\{ x^{k_{j}}\right\} _{j=0}^{\infty }$of $\left\{ x^{k}\right\}
_{k=0}^{\infty }$ which converges weakly to $x$.

For each point $x\in \mathcal{H},$\ there exists a unique nearest point in $%
D $, denoted by $P_{D}(x)$. That is,%
\begin{equation}
\left\Vert x-P_{D}\left( x\right) \right\Vert \leq \left\Vert x-y\right\Vert
\text{ for all }y\in D.
\end{equation}%
The mapping $P_{D}:\mathcal{H}\rightarrow D$ is called the metric projection
of $\mathcal{H}$ onto $D$. It is well known that $P_{D}$ is a \textit{%
nonexpansive mapping} of $\mathcal{H}$ onto $D$, and further more \textit{%
firmly nonexpansive mapping}. This is captured in the next lemma.

\begin{lemma}
\label{lem21} For any $x,y\in \mathcal{H}$ and $z\in D$, it holds

\begin{itemize}
\item $\Vert P_{D}(x)-P_{D}(y)\Vert\leq \Vert x-y\Vert ;\medskip $

\item $\Vert P_{D}(x)-z\Vert ^{2}\leq \Vert x-z\Vert ^{2}-\Vert
P_{D}(x)-x\Vert ^{2}$;
\end{itemize}
\end{lemma}

The characterization of the metric projection $P_{D}$ \cite[Section 3]%
{Goebel+Reich}, is given in the next lemma.

\begin{lemma}
\label{lem22} Let $x\in \mathcal{H}$ and $z\in D$. Then $z=P_{D}\left(
x\right) $ if and only if
\begin{equation}
P_{D}(x)\in D
\end{equation}%
and%
\begin{equation}
\left\langle x-P_{D}\left( x\right) ,P_{D}\left( x\right) -y\right\rangle
\geq 0\text{ for all }x\in \mathcal{H},\text{ }y\in D.
\end{equation}
\end{lemma}

\begin{definition}
The \texttt{normal cone} of $D$ at $v\in D$, denote by $N_{D}\left( v\right)
$ is defined as%
\begin{equation}
N_{D}\left( v\right) :=\{d\in \mathcal{H}\mid \left\langle
d,y-v\right\rangle \leq 0\text{ for all }y\in D\}.  \label{eq:normal_cone}
\end{equation}
\end{definition}

\begin{definition}
Let $B:\mathcal{H}\rightrightarrows 2^{\mathcal{H}}\mathcal{\ }$be a
point-to-set operator defined on a real Hilbert space $\mathcal{H}$. The
operator $B$ is called a \texttt{maximal monotone operator} if $B$ is
\texttt{monotone}, i.e.,%
\begin{equation}
\left\langle u-v,x-y\right\rangle \geq 0\text{ for all }u\in B(x)\text{ and }%
v\in B(y),
\end{equation}%
and the graph $G(B)$ of $B,$%
\begin{equation}
G(B):=\left\{ \left( x,u\right) \in \mathcal{H}\times \mathcal{H}\mid u\in
B(x)\right\} ,
\end{equation}%
is not properly contained in the graph of any other monotone operator.
\end{definition}

It is clear (\cite[Theorem 3]{Rockafellar}) that a monotone mapping $B$ is
maximal if and only if, for any $\left( x,u\right) \in\mathcal{H}\times%
\mathcal{H},$ if $\left\langle u-v,x-y\right\rangle \geq0$ for all $\left(
v,y\right) \in G(B)$, then it follows that $u\in B(x).$

\begin{lemma}
\cite{BC} \label{lem24} Let $D$ be a nonempty, closed and convex subset of
a Hilbert space $\mathcal{H}$. Let $\{x^{k}\}_{k=0}^{\infty }$ be a bounded
sequence which satisfies the following properties:

\begin{itemize}
\item every limit point of $\{x^{k}\}_{k=0}^{\infty }$ lies in $D$;

\item $\lim_{n\rightarrow \infty }\Vert x^{k}-x\Vert $ exists for every $%
x\in D$.
\end{itemize}
Then $\{x^{k}\}_{k=0}^{\infty }$ weakly converges to a point in $D$.
\end{lemma}

\begin{lemma}
\label{le4}
Assume that $\{a_k\}_{k=0}^{\infty }$ is a sequence of nonnegative real numbers such that
\begin{equation}
a_{k+1}\leq(1+\gamma_k)a_k+\delta_k,\,\,\,\,\,\,\forall k\geq0,
\end{equation}
where the sequences $\{\gamma_k\}_{k=0}^{\infty }\subset[0,+\infty)$ and $\{\delta_k\}_{k=0}^{\infty }$ satisfy

\begin{itemize}
\item
$\sum_{k=0}^\infty \gamma_k<+\infty$;

\item  $\sum_{k=0}^\infty \delta_k<+\infty$ or $\sup\delta_k\leq0$.
\end{itemize}
Then $\lim_{k\rightarrow\infty}a_k$ exists.
\end{lemma}

\textit{ Proof.} We prove the lemma only for $\sum_{k=0}^\infty \delta_k<+\infty$, when $\sup\delta_k\leq0$, the proof is similar.

For any natural number $l$ such that $1<l<k$, we have
\begin{equation}
\label{le4-1}
\aligned
a_{k+1}&\leq(1+\gamma_k)a_k+\delta_k\\
&\leq(1+\gamma_k)(1+\gamma_{k-1})a_{k-1}+(1+\gamma_k)\delta_{k-1}+\delta_k\\
&\leq(1+\gamma_k)(1+\gamma_{k-1})\cdots (1+\gamma_{l})a_{l}+(1+\gamma_k)\cdots(1+\gamma_{l+1})\delta_{l}\\
&\quad+\cdots+(1+\gamma_k)\delta_{k-1}+\delta_k\\
&\leq e^{\sum_{m=l}^{k}\gamma_m}\left(a_l+\sum_{m=l}^{k}\delta_m\right).
\endaligned
\end{equation}
Now fix $l$ and take superior limit for $k$:
\begin{equation}
\label{le4-2}
\varlimsup_{k\rightarrow\infty}a_{k}\leq e^{\sum_{m=l}^{+\infty}\gamma_m}\left(a_l+\sum_{m=l}^{+\infty}\delta_m\right).
\end{equation}
Thus,
\begin{equation}
\label{le4-3}
a_l\geq e^{-\sum_{m=l}^{+\infty}\gamma_m}\varlimsup_{k\rightarrow\infty}a_{k} -\sum_{m=l}^{+\infty}\delta_m.
\end{equation}
By taking now inferior limit for $l$ in the inequality (\ref{le4-3}) with $\sum_{k=0}^\infty \gamma_k<+\infty$ and $\sum_{k=0}^\infty \delta_k<+\infty$, we get
\begin{equation}
\label{le4-4}
\varliminf_{k\rightarrow\infty}a_k\geq \varlimsup_{k\rightarrow\infty}a_{k},
\end{equation}
which yields the existence of  $\lim_{k\rightarrow\infty}a_k$. $\Box$\\

Another useful property which derives easily from the Cauchy-Schwarz inequality and the mean value inequality is the following lemma.
\begin{lemma}
\label{lem-ineq}
Let $a,b\in \mathcal{H}$, then
\begin{equation}
\label{lem-ineq-1}
\left|2\langle a,b\rangle\right|\leq \|b\|\|a\|^2+\|b\|.
\end{equation}
\end{lemma}

\section{Convergence of the PC-algorithms with outer perturbations }\label{sec:Convergence}

In this section, we present two PC-algorithms with outer  perturbations and analyze their convergence. We first discuss the PC-algorithm I with outer perturbations.

\begin{algorithm} (PC-algorithm I with outer perturbations)\\
\label{PC1-OP}
Choose an arbitrary starting point $x^0\in \mathcal{H}.$
 Given the current iterate $x^k\in \mathcal{H}$, compute
\begin{equation}
\label{PC1-OP-y}
y^k=P_C(x^k-\beta_k F(x^k))+e_1(x^k),
\end{equation}
where $\beta_k>0$ is selected such that
\begin{equation}
\label{PC1-OP-beta}
\aligned
\beta_k\|F(x^k)-F(y^k-e_1^k)\|\leq\nu\|x^k-y^k+e_1^k\|,\,\,\,\, \nu\in(0,1).
\endaligned
\end{equation}
Define
\begin{equation}
\label{PC1-OP-d}
d(x^k,y^k)=(x^k-y^k+e_1(x^k))-\beta_k(F(x^k)-F(y^k-e_1(x^k))),
\end{equation}
and calculate
\begin{equation}
\label{PC1-OP-x}
x^{k+1}=x^k-\gamma\rho_k d(x^k,y^k)+e_2(x^k),
\end{equation}
where $\gamma\in(0,2)$,
and
\begin{equation}
\label{PC1-OP-rho}
\aligned
\rho_k:=\frac{\varphi(x^k,y^k)}{\|d(x^k,y^k)\|^2},
\endaligned
\end{equation}
where $\varphi(x^k,y^k)=\langle x^k-y^k+e_1(x^k), d(x^k,y^k)\rangle$.
\end{algorithm}

For the convergence proof we assume that the sequences of perturbations $\{e_i(x^k)\}_{k=0}^\infty$, $i=1,2$, are summable, i.e.,
\begin{equation}
\label{PC1-OP-Err}
\sum_{k=0}^\infty\|e_i(x^k)\|<+\infty,\quad i=1,2.
\end{equation}
For simplicity we denote $e_i^k:=e_i(x^k)$, $i=1,2$.

\begin{lemma}
\label{lem31}
Let $\{\rho_k\}_{k=0}^\infty$ be a sequence defined by (\ref{PC1-OP-rho}). Then under Conditions \ref{con:Condition 1.2} and \ref{con:Condition 1.3}, we have\\
\begin{equation}
\label{le11}
\rho_k\geq\frac{1-\nu}{1+\nu^2}.
\end{equation}
\end{lemma}

\textit{ Proof.}
From the Cauchy-Schwarz inequality and Condition \ref{con:Condition 1.3}, it follows
\begin{equation}
\label{le12}
\aligned
\varphi(x^k,y^k)&=\langle x^k-y^k+e_1^k,d(x^k,y^k)\rangle\\
&=\langle x^k-y^k+e_1^k,(x^k-y^k+e_1^k)-\beta_k(F(x^k)-F(y^k-e_1^k))\rangle\\
&=\|x^k-y^k+e_1^k\|^2-\beta_k\langle x^k-y^k+e_1^k,F(x^k)-F(y^k-e_1^k)\rangle\\
&\geq\|x^k-y^k+e_1^k\|^2-\beta_k\| x^k-y^k+e_1^k\|\|F(x^k)-F(y^k-e_1^k)\|\\
&\geq(1-\nu)\|x^k-y^k+e_1^k\|^2.\\
\endaligned
\end{equation}
Using Conditions \ref{con:Condition 1.2} and \ref{con:Condition 1.3}, we obtain
\begin{equation}
\label{le13}
\aligned
\|d(x^k,y^k)\|^2&=\|x^k-y^k+e_1^k-\beta_k(F(x^k)-F(y^k-e_1(x^k)))\|^2\\
&=\|x^k-y^k+e_1^k\|^2+\beta_k^2\|F(x^k)-F(y^k-e_1(x^k))\|^2\\
&\quad-2\beta_k\langle x^k-y^k+e_1^k,F(x^k)-F(y^k-e_1(x^k))\rangle\\
&\leq(1+\nu^2)\|x^k-y^k+e_1^k\|^2.\\
\endaligned
\end{equation}
Combining (\ref{le12}) and (\ref{le13}), we obtain (\ref{le11}) and the proof is complete. $\Box$

\begin{theorem}
\label{th31}
Assume that Conditions \ref{con:Condition 1.1}--\ref{con:Condition 1.3} hold.   Then
any sequence $\{x^k\}_{k=0}^\infty$   generated by  Algorithm \ref{PC1-OP} converges weakly to
a solution of the variational inequality problem $(\ref{vip})$.
\end{theorem}

\textit{ Proof.}
Let $x^*\in SOL(C,F)$. By the definition of $x^{k+1}$, we have
\begin{equation}
\label{th31-1}
\aligned
&\|x^{k+1}-x^*\|^2\\
&=\|x^k-\gamma\rho_kd(x^k,y^k)+e_2^k-x^*\|^2\\
&=\|x^k-x^*\|^2+\|\gamma\rho_kd(x^k,y^k)-e_2^k\|^2-2\langle x^k-x^*,\gamma\rho_kd(x^k,y^k)-e_2^k\rangle\\
&=\|x^k-x^*\|^2+\|x^{k+1}-x^k\|^2-2\gamma\rho_k\langle x^k-x^*, d(x^k,y^k)\rangle+2\langle x^k-x^*,e_2^k\rangle.\\
\endaligned
\end{equation}
It holds
\begin{equation}
\label{th31-2}
\langle x^k-x^*, d(x^k,y^k)\rangle=\langle x^k-y^k+e_1^k, d(x^k,y^k)\rangle+\langle y^k-e_1^k-x^*, d(x^k,y^k)\rangle.
\end{equation}
By the definition of $y^k$ and Lemma \ref{lem22}, we get
\begin{equation}
\label{th31-3}
\langle y^k-e_1^k-x^*, x^k-y^k+e_1^k-\beta_k F(x^k)\rangle\geq0.
\end{equation}
From Condition \ref{con:Condition 1.2}, it follows
\begin{equation}
\label{th31-4}
\langle y^k-e_1^k-x^*, \beta_k F(y^k-e_1^k)-\beta_k F(x^*)\rangle\geq0.
\end{equation}
Since $x^*\in SOL(C,F)$ and $y^k-e_1^k\in C$, we get from (\ref{vip})
\begin{equation}
\label{th31-5}
\langle y^k-e_1^k-x^*, \beta_k F(x^*)\rangle\geq0.
\end{equation}
Adding up (\ref{th31-3})-(\ref{th31-5}), we obtain
\begin{equation}
\label{th31-6}
\langle y^k-e_1^k-x^*, d(x^k,y^k)\rangle\geq0.
\end{equation}
From (\ref{th31-2}), we get
\begin{equation}
\label{th31-7}
\langle x^k-x^*, d(x^k,y^k)\rangle\geq\langle x^k-y^k+e_1^k, d(x^k,y^k)\rangle=\varphi(x^k,y^k).
\end{equation}
Substituting (\ref{th31-7}) into (\ref{th31-1}), we get
\begin{equation}
\label{th31-8}
\aligned
\|x^{k+1}-x^*\|^2&\leq\|x^k-x^*\|^2+\|x^{k+1}-x^k\|^2-2\gamma\rho_k\varphi(x^k,y^k)\\
&\quad+2\langle x^k-x^*,e_2^k\rangle.\\
\endaligned
\end{equation}
By Lemma \ref{lem-ineq},
\begin{equation}
\label{th31-9}
\aligned
2\langle x^k-x^*,e_2^k\rangle\leq\|e_2^k\|+\|e_2^k\|\|x^k-x^*\|^2.\\
\endaligned
\end{equation}
Again using the definition of $x^{k+1}$, we have
\begin{equation}
\label{th31-10}
\aligned
-2\gamma\rho_k\varphi(x^k,y^k)&=-2\frac{1}{\gamma}\|\gamma\rho_kd(x^k,y^k)\|^2\\
&=-2\frac{1}{\gamma}\|x^k-x^{k+1}+e_2^k\|^2\\
&\leq-2\frac{1}{\gamma}\|x^k-x^{k+1}\|^2-4\frac{1}{\gamma}\langle x^k-x^{k+1},e_2^k\rangle\\
&\leq-2\frac{1}{\gamma}(1-\|e_2^k\|)\|x^k-x^{k+1}\|^2+2\frac{1}{\gamma}\|e_2^k\|,\\
\endaligned
\end{equation}
where the last inequality follows by Lemma \ref{lem-ineq}.
Adding (\ref{th31-8})-(\ref{th31-10}), we obtain
\begin{equation}
\label{th31-11}
\aligned
\|x^{k+1}-x^*\|^2&\leq(1+\|e_2^k\|)\|x^k-x^*\|^2-\frac{1}{\gamma}(2-\gamma-2\|e_2^k\|)\|x^{k+1}-x^k\|^2\\
&\quad+\frac{2+\gamma}{\gamma}\|e_2^k\|.
\endaligned
\end{equation}
From (\ref{PC1-OP-Err}), it follows
\begin{equation}
\label{th31-12}
\lim_{k\rightarrow\infty}\|e_i^k\|=0,\quad i=1,2.
\end{equation}
Therefore, we assume $\|e_2^k\|\in[0,1-\frac{\gamma}{2}-\mu)$, $k\geq0$, where $\mu\in(0,1-\frac{\gamma}{2})$. So, we get
\begin{equation}
\label{th31-13}
\aligned
\|x^{k+1}-x^*\|^2&\leq(1+\|e_2^k\|)\|x^k-x^*\|^2+\frac{2+\gamma}{\gamma}\|e_2^k\|-\frac{2\mu}{\gamma}\|x^{k+1}-x^k\|^2\\
&\leq(1+\|e_2^k\|)\|x^k-x^*\|^2+\frac{2+\gamma}{\gamma}\|e_2^k\|.\\
\endaligned
\end{equation}
Using (\ref{PC1-OP-Err}) and Lemma  \ref{le4}, the existence of the limit $\lim_{k\rightarrow\infty}\|x^k-x^*\|^2$ is guarantied and hence also the boundedness of the sequence $\{x^k\}_{k=0}^\infty$.

From (\ref{th31-13}) and the existence of $\lim_{k\rightarrow\infty}\|x^k-x^*\|^2$, it follows
\begin{equation}
\label{th31-14}
\sum_{k=0}^{\infty}\|x^{k+1}-x^k\|\leq+\infty
\end{equation}
which implies
\begin{equation}
\label{th31-15}
\lim_{k\rightarrow\infty}\|x^{k+1}-x^k\|=0.
\end{equation}
From (\ref{PC1-OP-x}), (\ref{PC1-OP-rho}) and Lemma \ref{lem31}, we have
\begin{equation}
\label{th31-16}
\aligned
\varphi(x^k,y^k)&=\frac{1}{\rho_k \gamma^2}\|x^k-x^{k+1}+e_2^k\|^2\\
&\leq\frac{1+\nu^2}{(1-\nu)\gamma^2}\|x^k-x^{k+1}+e_2^k\|^2\\
&\leq\frac{2(1+\nu^2)}{(1-\nu)\gamma^2}[\|x^k-x^{k+1}\|^2+\|e_2^k\|^2].\\
\endaligned
\end{equation}
Combining (\ref{le12}) and (\ref{th31-16}), we get
\begin{equation}
\label{th31-17}
\|x^k-y^k+e_1^k\|^2\leq\frac{2(1+\nu^2)}{(1-\nu)^2\gamma^2}\left[\|x^k-x^{k+1}\|^2+\|e_2^k\|^2\right].
\end{equation}
Using (\ref{th31-12}) and (\ref{th31-15}), we have
\begin{equation}
\label{th31-18}
\lim_{k\rightarrow\infty}\|x^k-y^k+e_1^k\|=0.
\end{equation}

Now, we show $\omega_w(x^k)\subseteq SOL(C,F).$
Due to the boundedness of the sequence $\{x^k\}_{k=0}^\infty$, it has at least one weak
accumulation point, we denote it by $\hat x\in \omega_w(x^k)$. So, there exists a subsequence $\{x^{k_i}\}_{i=0}^\infty$ of $\{x^k\}_{k=0}^\infty$
which converges weakly to $\hat x$. From (\ref{th31-18}), it follows that  $\{y^{k_i}-e_1^{k_i}\}_{i=0}^\infty$ also converges weakly to $\hat x.$
It is now left to show that $\hat x$ also solves the variational inequality (\ref{vip}).
Define the operator
\begin{equation}
\label{th31-19}
Av=
\left\{
\aligned
&f(v)+N_C(v), &v\in C,\\
&\emptyset, & v\notin C.
\endaligned
\right.
\end{equation}

It is known that $A$ is a maximal monotone operator and $A^{-1}(0) = SOL(C,F)$. If $(v,w) \in G(A)$,
 then we have $w-F(v) \in N_C(v)$ since $w\in A(v) = F(v) + N_C(v)$. Thus it follows that
\begin{equation}
\label{th31-21}
\langle w-F(v),v-y\rangle\geq0,\quad y\in C.
\end{equation}
Since $y^{k_i}-e_1^{k_i}\in C$, we have
\begin{equation}
\label{th31-22}
\langle w-F(v),v-y^{k_i}+e_1^{k_i}\rangle\geq0.
\end{equation}
On the other hand, by the definition of $y^k$ and Lemma \ref{lem22}, it follows that
\begin{equation}
\label{th31-23}
\langle x^k-\beta_k F(x^k)+e_1^k-y^k, y^{k}-e_1^k-v\rangle\geq0,
\end{equation}
and consequently,
\begin{equation}
\label{th31-24}
\left\langle \frac{y^k-e_1^k-x^k}{\beta_k}+F(x^k),v-y^k+e_1^k\right\rangle\geq0.
\end{equation}
Hence we have
\begin{equation}
\label{th31-25}
\aligned
&\langle w,v-y^{k_i}+e_1^{k_i}\rangle\\
&\geq\langle F(v),v-y^{k_i}+e_1^{k_i}\rangle\\
&\geq\langle F(v),v-y^{k_i}+e_1^{k_i}\rangle-\Big\langle \frac{y^{k_i}-e_1^{k_i}-x^{k_i}}{\beta_{k_i}}+F(x^{k_i}),v-y^{k_i}+e_1^{k_i}\Big\rangle\\
&=\langle F(v)-F(y^{k_i}-e_1^{k_i}),v-y^{k_i}+e_1^{k_i}\rangle+\langle F(y^{k_i}-e_1^{k_i})-F(x^{k_i}),v-y^{k_i}+e_1^{k_i}\rangle\\
&\quad -\Big\langle \frac{y^{k_i}-e_1^{k_i}-x^{k_i}}{\beta_{k_i}},v-y^{k_i}+e_1^{k_i}\Big\rangle\\
&\geq\langle F(y^{k_i}-e_1^{k_i})-F(x^{k_i}),v-y^{k_i}+e_1^{k_i}\rangle-\Big\langle \frac{y^{k_i}-e_1^{k_i}-x^{k_i}}{\beta_{k_i}},v-y^{k_i}+e_1^{k_i}\Big\rangle,\\
\endaligned
\end{equation}
which implies
\begin{equation}
\label{th31-26}
\aligned
\langle w,v-y^{k_i}+e_1^{k_i}\rangle&\geq\langle F(y^{k_i}-e_1^{k_i})-F(x^{k_i}),v-y^{k_i}+e_1^{k_i}\rangle\\
&\quad-\Big\langle \frac{y^{k_i}-e_1^{k_i}-x^{k_i}}{\beta_{k_i}},v-y^{k_i}+e_1^{k_i}\Big\rangle.\\
\endaligned
\end{equation}
Taking the limit as $i\rightarrow\infty$ in the above inequality, we obtain
\begin{equation}
\label{th31-27}
\langle w,v-\hat x\rangle\geq0.
\end{equation}
Since $A$ is a maximal monotone operator, it follows that $\hat x\in A^{-1}(0) =
SOL(C,F)$. So, $\omega_w(x^k)\subseteq SOL(C,F).$
Finally, since $\lim_{k\rightarrow\infty}\|x^{k}-x^*\|$ exists, $\omega_w(x^k)\subseteq SOL(C,F)$ and by using Lemma \ref{lem24}, we conclude that $\{x^k\}_{k=0}^\infty$ weakly converges to a solution of the variational inequality (\ref{vip}), which completes the proof. $\Box$\\

Now that we proved the converges of the PC-algorithm I with outer perturbations, we follow Cai {\it et al.} \cite{He1} and show that that it converges at a $O(1/t)$ rate.

\begin{lemma}
\label{lem35}
Let $\{x^k\}_{k=0}^\infty$ and $\{y^k\}_{k=0}^\infty$ be any two sequences generated by Algorithm \ref{PC1-OP}. Then we have
\begin{equation}
\label{lem35-0}
\aligned
\langle x-y^k+e_1^k,\gamma\rho_k\beta_kF(y^k-e_1^k)\rangle&+
\frac12\left(\|x-x^k\|^2-\|x-x^{k+1}+e_2^k\|^2\right)\\
&\geq\frac12\gamma(2-\gamma)\rho_k^2\|d(x^k,y^k)\|^2,\quad \forall x\in C.
\endaligned
\end{equation}
\end{lemma}

\textit{Proof.}
Notice that the projection equation (\ref{PC1-OP-y}) can be written as
\begin{equation}
\label{lem35-1}
y^k-e_1^k=P_C(y^k-e^k_1-(\beta_kF(y^k-e_1^k)-d(x^k,y^k))).
\end{equation}
From Lemma \ref{lem22} we have
\begin{equation}
\label{lem35-2}
\langle x- y^k+e^k_1,\beta_kF(y^k-e_1^k)-d(x^k,y^k)\rangle\geq0,\quad \forall x\in C,
\end{equation}
which implies
\begin{equation}
\label{lem35-2a}
\langle x- y^k+e^k_1,\beta_kF(y^k-e_1^k)\rangle
\geq\langle x- y^k+e^k_1,d(x^k,y^k)\rangle,\quad \forall x\in C.
\end{equation}
Due to (\ref{PC1-OP-x}), we have
\begin{equation}
\label{lem35-4}
\gamma\rho_k d(x^k,y^k)=x^k-x^{k+1}+e_2^k,
\end{equation}
which with (\ref{lem35-2a}) yields
\begin{equation}
\label{lem35-5}
\aligned
\langle x- y^k+e^k_1,\gamma\rho_k\beta_kF(y^k-e_1^k)\rangle
&\geq\langle x- y^k+e^k_1,x^k-x^{k+1}+e_2^k\rangle.
\endaligned
\end{equation}
Now using the following identity for (\ref{lem35-5})
\begin{equation}
\langle a-b,c-d\rangle=\frac12(\|a-d\|^2-\|a-c\|^2)+\frac12(\|c-b\|^2-\|d-b\|^2),
\end{equation}
we obtain
\begin{equation}
\label{lem35-6}
\aligned
\langle x- (y^k-e^k_1),&x^k-(x^{k+1}-e_2^k)\rangle\\
&=\frac12(\|x-x^{k+1}+e_2^k\|^2-\|x-x^k\|^2)\\
&\quad+\frac12(\|x^k+e_1^k-y^k\|^2-\|x^{k+1}-y^k-e_2^k+e_1^k\|^2).
\endaligned
\end{equation}
By using $x^{k+1}=x^k-\gamma\rho_k d(x^k,y^k)+e_2(x^k),$  and $(\ref{PC1-OP-rho})$, we get
\begin{equation}
\label{lem35-7}
\aligned
\|x^k+e_1^k&-y^k\|^2-\|x^{k+1}-y^k-e_2^k+e_1^k\|^2\\
&=\|x^k+e_1^k-y^k\|^2-\|(x^{k}-y^k+e_1^k)-\gamma\rho_k d(x^k,y^k)\|^2\\
&=2\gamma\rho_k \langle x^k+e_1^k-y^k, d(x^k,y^k)\rangle-\gamma^2\rho_k^2\| d(x^k,y^k)\|^2\\
&=\gamma(2-\gamma)\rho_k^2\| d(x^k,y^k)\|^2.
\endaligned
\end{equation}
Combining (\ref{lem35-5}), (\ref{lem35-6}) and (\ref{lem35-7}), we get (\ref{lem35-0}),
and the desired result is obtained. $\Box$

\begin{theorem}
\label{Th32} Assume that Conditions \ref{con:Condition 1.1}--\ref{con:Condition 1.3} hold. Let $%
\{x^{k}\}_{k=0}^{\infty }$ and $\{y^{k}\}_{k=0}^{\infty }$ be any two sequences generated by
Algorithm \ref{PC1-OP}. For any integer $t>0$, there exists a point $\hat y_{t}\in C$
such that
\begin{equation}
\label{Th32-0a}
\aligned
\langle F(x),\hat y_{t}-x\rangle \leq \frac{1}{2\gamma \Upsilon _{t}}%
\left(\Vert x-x^{0}\Vert ^{2}+2M\right),\quad \forall x\in C,
\endaligned
\end{equation}%
where
\begin{equation}
\label{Th32-0b}
\aligned
\hat y_{t}=\frac{1}{\Upsilon _{t}}\sum_{k=0}^{t}\rho_k\beta_k(y^{k}-e^k_1),\,\, \Upsilon
_{t}=\sum_{k=0}^{t}\rho_k\beta_k,
 \hbox{ and }\,\, M=\sup_{k\in {\Bbb N}}\{\Vert x^{k+1}-x\Vert
\}\sum_{k=0}^{\infty }\Vert e_{2}^{k}\Vert.
\endaligned
\end{equation}
Further, we also have
\begin{equation}
\label{Th32-0c}
\aligned\langle F(x),y_{t}-x\rangle \leq \frac{1}{2\gamma \Upsilon _{t}}%
(\Vert x-x^{0}\Vert ^{2}+2M)
+\frac{\|F(x)\|}{\Upsilon _{t}}\sum_{k=0}^t\rho_k\beta_k\|e_1^k\|,\quad \forall x\in C,\endaligned
\end{equation}%
where $\Upsilon _{t}$ and $M$ are defined as in (\ref{Th32-0b}), and
\begin{equation}
\label{Th32-0d}
y_{t}=\frac{1}{\Upsilon _{t}}\sum_{k=0}^{t}\rho_k\beta_ky^{k}.
\end{equation}
\end{theorem}

\textit{Proof.} Take an arbitrary point $x\in C.$ By Condition \ref{con:Condition 1.2}, we have
\begin{equation}
\label{Th32-1}
\langle x-y^k+e_1^k,\rho_k\beta_kF(y^k-e_1^k)\rangle
\leq\langle  x-y^k+e_1^k, \rho_k\beta_kF(x)\rangle,
\end{equation}
which with (\ref{lem35-0}) implies
\begin{equation}
\label{Th32-2}
\aligned
\langle  y^k-e_1^k&-x, \rho_k\beta_kF(x)\rangle\\
&\leq
\frac1{2\gamma}\left(\|x-x^k\|^2-\|x-x^{k+1}+e_2^k\|^2\right)\\
&\leq
\frac1{2\gamma}\left(\|x-x^k\|^2-\|x-x^{k+1}\|^2-2\langle x-x^{k+1},e_2^k\rangle \right)\\
&\leq
\frac1{2\gamma}\left(\|x-x^k\|^2-\|x-x^{k+1}\|^2+2\| x-x^{k+1}\|\|e_2^k\| \right).\\
\endaligned
\end{equation}
Summing the inequalities (\ref{Th32-2}) over $k = 0,\ldots,t$, we obtain
\begin{equation}
\label{Th32-3}
\aligned \left\langle
\sum_{k=0}^t\rho_k\beta_k(y^k-e^k_1)-\left(\sum_{k=0}^t\rho_k\beta_k\right)x,F(x)\right\rangle
&\leq\frac1{2\gamma}\|x^0-x\|^2+\frac{M_1}\gamma\sum_{k=0}^t\|e^k_2\|, \\
\endaligned
\end{equation}
where $M_1=\sup_{k\in {\Bbb N}}\| x-x^{k+1}\|.$
Using the notations of $\Upsilon_t$ and $\hat y_t$ in the above inequality, we
derive
\begin{equation}
\label{Th32-4}
\langle F(x),\hat y_t-x\rangle\leq\frac{1}{2\gamma\Upsilon_t}(\|x-x^0\|^2+2M),%
\quad \forall x\in C.
\end{equation}
From (\ref{Th32-3}), it follows
\begin{equation}
\label{Th32-5}
\aligned
&\left\langle
\sum_{k=0}^t\rho_k\beta_ky^k-\left(\sum_{k=0}^t\rho_k\beta_k\right)x,F(x)\right\rangle\\
&\leq\frac1{2\gamma}\|x^0-x\|^2+\frac{M_1}\gamma\sum_{k=0}^t\|e^k_2\|
+ \left\langle\sum_{k=0}^t\rho_k\beta_ke_1^k,F(x)\right\rangle, \\
&\leq\frac1{2\gamma}\|x^0-x\|^2+\frac{M_1}\gamma\sum_{k=0}^t\|e^k_2\|
+ \|F(x)\|\sum_{k=0}^t\rho_k\beta_k\|e_1^k\|. \\
\endaligned
\end{equation}
Similarly with (\ref{Th32-4}), we get (\ref{Th32-0c})
and the desired result is obtained. $\Box$

\begin{remark}
\textrm{\ From Lemma \ref{lem31}, it follows that
\begin{equation}  \label{m5}
\Upsilon_t\geq(t+1)\underline{\gamma}.
\end{equation}
So, due to (\ref{Th32-0a}), we get that Algorithm \ref{PC1-OP} converges at the rate of $O(1/t)$.}
\end{remark}

Next we wish to study the convergence (also its rate) of the PC-algorithm II with outer perturbations. The analysis follows similar lines as the one presented earlier for the PC-algorithm I, but it is presented next in full details for the convenience of the reader.
\begin{algorithm} (PC-algorithm II with outer perturbations)\\
\label{PC2-OP}
Choose an arbitrary starting point $x^0\in \mathcal{H}.$
 Given the current iterate $x^k\in \mathcal{H}$, compute
\begin{equation}
\label{PC2-OP-y}
y^k=P_C(x^k-\beta_k F(x^k)+e_1(x^k)),
\end{equation}
where
$\beta_k>0$ is selected such that
\begin{equation}
\label{PC2-OP-y-beta}
\aligned
\beta_k\|F(x^k)-F(y^k)\|\leq\nu\|x^k-y^k\|,\,\,\,\, \nu\in(0,1).
\endaligned
\end{equation}
Caculate
\begin{equation}
\label{PC2-OP-x}
x^{k+1}=P_C(x^k-\gamma\rho_k\beta_k F(y^k)+e_2(x^k)),
\end{equation}
where $\gamma\in(0,2)$,
\begin{equation}
\label{PC2-OP-rho}
\aligned
\rho_k:=\frac{\langle x^k-y^k, d(x^k,y^k)\rangle}{\|d(x^k,y^k)\|^2},
\endaligned
\end{equation}
and
\begin{equation}
\label{PC2-OP-d}
d(x^k,y^k):=(x^k-y^k)-\beta_k(F(x^k)-F(y^k))+e_1(x^k).
\end{equation}
\end{algorithm}

As previously, we assume that $e_1(x^k)$ and $e_2(x^k)$ satisfy (\ref{PC1-OP-Err}), and in addition we also need to assume that
\begin{equation}
\label{PC2-OP-e1}
\|e_1(x^k)\|\leq\mu\|x^k-y^k\|,
\end{equation}
where $\mu\in[0,1-\nu)$.

\begin{lemma}
\label{lem32}
Let $\{\rho_k\}_{k=0}^\infty$ be a sequence which is defined by (\ref{PC2-OP-rho}). Then under Conditions \ref{con:Condition 1.2} and \ref{con:Condition 1.3}, we have
\begin{equation}
\label{lem32-rho}
\rho_k\geq \frac{1-\nu-\mu}{1+\nu^2+\mu^2+2\mu+2\nu\mu}.
\end{equation}
\end{lemma}

\textit{ Proof.}
By the definition of $d(x^k,y^k)$, we get
\begin{equation}
\label{lem32-1}
\aligned
&\langle x^k-y^k, d(x^k,y^k)\rangle\\
&=\|x^k-y^k\|^2-\beta_k\langle x^k-y^k,F(x^k)-F(y^k)\rangle+\langle e_1^k,x^k-y^k\rangle\\
&\geq\|x^k-y^k\|^2-\beta_k\| x^k-y^k\|\|F(x^k)-F(y^k)\|-\| e_1^k\|\|x^k-y^k\|\\
&\geq\|x^k-y^k\|^2-\nu\| x^k-y^k\|^2-\mu\|x^k-y^k\|^2\\
&\geq(1-\nu-\mu)\|x^k-y^k\|^2.\\
\endaligned
\end{equation}
On the other hand,
\begin{equation}
\label{lem32-2}
\aligned
&\| d(x^k,y^k)\|^2\\
&=\|x^k-y^k\|^2+\beta_k^2\|F(x^k)-F(y^k)\|^2+\|e_1^k\|^2+2\langle e_1^k,x^k-y^k\rangle\\
&\quad -2\beta_k\langle x^k-y^k,F(x^k)-F(y^k)\rangle-2\beta_k\langle e_1^k,F(x^k)-F(y^k)\rangle\\
&\leq(1+\nu^2+\mu^2)\|x^k-y^k\|^2+2\|e_1^k\|\| x^k-y^k\|+2\beta_k\|e_1^k\|\|F(x^k)-F(y^k)\|\\
&\leq(1+\nu^2+\mu^2)\|x^k-y^k\|^2+2\mu\| x^k-y^k\|^2+2\nu\mu\|x^k-y^k\|^2\\
&\leq(1+\nu^2+\mu^2+2\mu+2\nu\mu)\|x^k-y^k\|^2.\\
\endaligned
\end{equation}
So, we get (\ref{lem32-rho}), and the proof is complete. $\Box$

\begin{theorem}
\label{th33}
Assume that Conditions \ref{con:Condition 1.1}--\ref{con:Condition 1.3} hold. Then
any sequence $\{x^k\}_{k=0}^\infty$ generated by  Algorithm \ref{PC2-OP} converges weakly to
a solution of the variational inequality problem $(\ref{vip})$.
\end{theorem}

\textit{ Proof.}
Let $x^*\in SOL(C,F)$. By the definition of $x^{k+1}$ and Lemma \ref{lem21}, we have
\begin{equation}
\label{th32-1}
\aligned
&\|x^{k+1}-x^*\|^2\\
&\leq\|x^k-\gamma\rho_k\beta_kF(y^k)+e_2^k-x^*\|^2-\|x^k
-\gamma\rho_k\beta_kF(y^k)+e_2^k-x^{k+1}\|^2\\
&=\|x^k-x^*\|^2-\|x^{k+1}-x^k\|^2-2\langle x^{k+1}-x^*,\gamma\rho_k\beta_kF(y^k)-e_2^k\rangle\\
&=\|x^k-x^*\|^2-\|x^{k+1}-x^k\|^2-2\gamma\rho_k\beta_k\langle x^{k+1}-x^*, F(y^k)\rangle\\
&\quad+2\langle x^{k+1}-x^*,e_2^k\rangle.\\
\endaligned
\end{equation}
Notice that the projection equation (\ref{PC2-OP-y}) can be written as
\begin{equation}
\label{th32-1a}
y^k=P_C(y^k-(\beta_kF(y^k)-d(x^k,y^k))).
\end{equation}
So, from Lemma \ref{lem22} we have
\begin{equation}
\label{th32-1b}
\langle x- y^k,\beta_kF(y^k)-d(x^k,y^k)\rangle\geq0,\quad \forall x\in C,
\end{equation}
which with $x^{k+1}\in C$ implies
\begin{equation}
\label{th32-1c}
\langle x^{k+1}- y^k,\beta_kF(y^k)-d(x^k,y^k)\rangle\geq0.
\end{equation}
Since $x^*\in SOL(C,F)$ and $y^k\in C$ we get from (\ref{vip})
\begin{equation}
\label{th32-1d}
\langle y^k-x^*, \beta_k F(x^*)\rangle\geq0.
\end{equation}
Using (\ref{th32-1c}) and (\ref{th32-1d}), we get
\begin{equation}
\label{th32-2}
\aligned
&-2\gamma\rho_k\beta_k\langle x^{k+1}-x^*, F(y^k)\rangle\\
&=-2\gamma\rho_k\beta_k\langle x^{k+1}-y^k, F(y^k)\rangle-2\gamma\rho_k\beta_k\langle y^k-x^*, F(y^k)\rangle\\
&\leq-2\gamma\rho_k\langle x^{k+1}-y^k, d(x^k,y^k)\rangle\\
&=-2\gamma\rho_k\langle x^k-y^k, d(x^k,y^k)\rangle+2\gamma\rho_k\langle x^k-x^{k+1}, d(x^k,y^k)\rangle\\
&\leq-2\gamma\rho_k^2\|d(x^k,y^k)\|^2+2\gamma\rho_k\|x^k-x^{k+1}\|\|d(x^k,y^k)\|\\
&\leq-\gamma(2-\gamma)\rho_k^2\|d(x^k,y^k)\|^2+\|x^k-x^{k+1}\|^2.\\
\endaligned
\end{equation}
By Lemma \ref{lem-ineq},
\begin{equation}
\label{th32-3}
2\langle x^{k+1}-x^*,e_2^k\rangle\leq\|e_2^k\|+\|e_2^k\|\|x^{k+1}-x^*\|^2.
\end{equation}
Adding (\ref{th32-1}), (\ref{th32-2}) and (\ref{th32-3}), we obtain
\begin{equation}
\label{th32-4}
\aligned
(1-\|e_2^k\|)\|x^{k+1}-x^*\|^2&\leq\|x^k-x^*\|^2-\gamma(2-\gamma)\rho_k^2\|d(x^k,y^k)\|^2+\|e_2^k\|.
\endaligned
\end{equation}
From (\ref{PC1-OP-Err}), it follows
\begin{equation}
\label{th32-5}
\lim_{k\rightarrow\infty}\|e_i^k\|=0,\quad i=1,2.
\end{equation}
Therefore, we assume $\|e_2^k\|\in[0,1/2)$, $k\geq0$. So,
\begin{equation}
\label{th32-7}
1\leq\frac{1}{1-\|e_2^k\|}\leq1+2\|e_2^k\|\leq2.
\end{equation}
By (\ref{th32-4}) and (\ref{th32-7}) we have
\begin{equation}
\label{th32-8}
\aligned
&\|x^{k+1}-x^*\|^2\\
&\leq\frac{1}{1-\|e_2^k\|}\|x^k-x^*\|^2-
\frac{\gamma(2-\gamma)\rho_k^2}{1-\|e_2^k\|}\|d(x^k,y^k)\|^2+\frac{\|e_2^k\|}{1-\|e_2^k\|}\\
&\leq(1+2\|e_2^k\|)\|x^k-x^*\|^2-\gamma(2-\gamma)
\rho_k^2\|d(x^k,y^k)\|^2+\|e_2^k\|(1+2\|e_2^k\|)\\
&\leq(1+2\|e_2^k\|)\|x^k-x^*\|^2+2\|e_2^k\|.\\
\endaligned
\end{equation}
Following the proof of (\ref{th31-15}), we get
\begin{equation}
\label{th32-10}
\lim_{k\rightarrow\infty}\rho_k^2\|d(x^k,y^k)\|^2=0.
\end{equation}
From (\ref{lem32-1}) and Lemma \ref{lem32}, we get
\begin{equation}
\label{th32-12}
\aligned
(1-\nu-\mu)\| x^k-y^k\|^2&\leq\langle x^k-y^k,d(x^k,y^k)\rangle\\
&=\rho_k\|d(x^k,y^k)\|^2\\
&\leq\frac{1+\nu^2+\mu^2+2\mu+2\nu\mu}{1-\nu-\mu}\rho_k^2\|d(x^k,y^k)\|^2,\\
\endaligned
\end{equation}
which with (\ref{th32-10}) yields
\begin{equation}
\label{th32-13}
\lim_{k\rightarrow\infty}\|x^k-y^k\|^2=0.
\end{equation}
Now the rest of the proof follows directly the proof of Theorem \ref{th31}, and therefore we obtain the desired result. $\Box$\\

The next step is to evaluate the convergence rate of Algorithm \ref{PC2-OP}.

\begin{lemma}
\label{lem311}
Let $\{x^k\}_{k=0}^\infty$ and $\{y^k\}_{k=0}^\infty$ be given by Algorithm \ref{PC2-OP}. Then we have
\begin{equation}
\label{lem311-0}
\aligned
&\langle x-y^k,\gamma\rho_k\beta_kF(y^k)\rangle+
\frac12\left(\|x-x^k\|^2-\|x-x^{k+1}\|^2\right)\\
&\geq\frac12\gamma(2-\gamma)\rho_k^2\|d(x^k,y^k)\|^2+
\langle x-x^{k+1},e_2^k\rangle,\quad \forall x\in C.
\endaligned
\end{equation}
\end{lemma}

\textit{ Proof.}
Using (\ref{th32-1c}), we get
\begin{equation}
\label{lem311-2}
\aligned
\langle x^{k+1}- y^k,\gamma\rho_k\beta_kF(y^k)\rangle
&\geq\gamma\rho_k\langle x^{k+1}- y^k,d(x^k,y^k)\rangle\\
&=\gamma\rho_k\langle x^{k}- y^k,d(x^k,y^k)\rangle\\
&\quad-\gamma\rho_k\langle x^{k}- x^{k+1},d(x^k,y^k)\rangle.\\
\endaligned
\end{equation}
In order to evaluate the last term of (\ref{lem311-2}), we use (\ref{PC2-OP-rho}) and get
\begin{equation}
\label{lem311-3}
\gamma\rho_k\langle x^{k}- y^k,d(x^k,y^k)\rangle
=\gamma\rho_k^2\|d(x^k,y^k)\|^2.
\end{equation}
Using the Cauchy-Schwarz inequality we get
\begin{equation}
\label{lem311-4}
-\gamma\rho_k\langle x^{k}- x^{k+1},d(x^k,y^k)\rangle\geq-\frac12\|x^{k}- x^{k+1}\|^2
-\frac12\gamma^2\rho_k^2\|d(x^k,y^k)\|^2.
\end{equation}
and obtain
\begin{equation}
\label{lem311-4a}
\langle x^{k+1}- y^k,\gamma\rho_k\beta_kF(y^k)\rangle\geq
\frac12\gamma(2-\gamma)\rho_k\|d(x^k,y^k)\|^2-\frac12\|x^{k}- x^{k+1}\|^2.
\end{equation}
By Lemma \ref{lem22} and (\ref{PC2-OP-x}), we have
\begin{equation}
\label{lem311-5}
\langle x^k-\gamma\rho_k\beta_k F(y^k)+e_2^k-x^{k+1},x-x^{k+1}\rangle\leq0,\quad \forall x\in C,
\end{equation}
and consequently
\begin{equation}
\label{lem311-6}
\gamma\rho_k\beta_k\langle  F(y^k),x-x^{k+1}\rangle\geq
\langle x^k+e_2^k-x^{k+1},x-x^{k+1}\rangle,\quad \forall x\in C.
\end{equation}
Using the identity $\langle a,b\rangle=\frac12(\|a\|^2-\|a-b\|^2+\|b\|^2)$
for the right hand side of (\ref{lem311-6}), we obtain
\begin{equation}
\label{lem311-7}
\aligned
\gamma\rho_k\beta_k&\langle  F(y^k),x-x^{k+1}\rangle\\
&\geq\frac12(\|x^k+e_2^k-x^{k+1}\|^2-\|x^k+e_2^k-x\|^2+\|x-x^{k+1}\|^2)\\
&=\frac12(\|x^k-x^{k+1}\|^2-\|x^k-x\|^2+\|x-x^{k+1}\|^2)
+\langle x-x^{k+1},e_2^k\rangle.\\
\endaligned
\end{equation}
Adding (\ref{lem311-4a}) and (\ref{lem311-7}), we get (\ref{lem311-0}) and the proof
is complete. $\Box$\\

Now, in the same spirit of Theorem \ref{Th32}, by using Lemma \ref{lem311}, the convergence rate (O($1/$t)) of
Algorithm \ref{PC2-OP} is guaranteed.

\begin{theorem}
\label{Th312} Assume that Conditions \ref{con:Condition 1.1}--\ref{con:Condition 1.3} hold. Let $%
\{x^{k}\}_{k=0}^{\infty }$ and $\{y^{k}\}_{k=0}^{\infty }$ be any sequences generated by
Algorithm \ref{PC1-OP}. For any integer $t>0$, we have a $y_{t}\in C$
which satisfies
\begin{equation}
\label{Th312-0a}
\aligned
\langle F(x), y_{t}-x\rangle \leq \frac{1}{2\gamma \Upsilon _{t}}%
(\Vert x-x^{0}\Vert ^{2}+2M),\quad \forall x\in C,\endaligned
\end{equation}%
where
\begin{equation}
\label{Th312-0b}
y_{t}=\frac{1}{\Upsilon _{t}}\sum_{k=0}^{t}\rho_k\beta_ky^{k},\,\, \Upsilon
_{t}=\sum_{k=0}^{t}\rho_k\beta_k,\,\, \hbox{and}\,\, M=\sup_{k\in {\Bbb N}}\{\Vert x^{k+1}-x\Vert
\}\sum_{k=0}^{\infty }\Vert e_{2}^{k}\Vert.
\end{equation}

\end{theorem}

\section{The bounded perturbation resilience of the PC-algorithms}\label{sec:BPR}

\qquad
In this section, we prove the bounded perturbation resilience (BPR) of the PC-algorithms.
This property is fundamental for the application of the superiorization methodology
(SM).

\subsection{Bounded perturbation resilience}

\qquad
The \textit{superiorization methodology} first appeared in Butnariu et al. in \cite{bdhk07}, without mentioning specifically the words superiorization and perturbation resilience. Some of the results in \cite{bdhk07} are based on earlier results of Butnariu, Reich and Zaslavski \cite{brz06, brz07, brz08}. For the state of current research on superiorization, visit the webpage:
\textquotedblleft Superiorization and Perturbation Resilience of Algorithms: A Bibliography compiled and continuously updated by Yair Censor\textquotedblright\  at: \href{http://math.haifa.ac.il/yair/bib-superiorization-censor.html}{http://math.haifa.ac.il/yair/bib-superiorization-censor.html}
and in particular see \cite[Section 3]{Censor14} and \cite[Appendix]{Censor17}.

Originally, the superiorization methodology is intended for constrained minimization (CM) problems of the form:
\begin{equation}
\label{e1}
\hbox{ minimize}\,\{\phi(x)\,|\, x\in \Psi\}
\end{equation}
where $\phi : \mathcal{H} \rightarrow \mathbb{R}$ is an objective function and $\Psi\subseteq \mathcal{H}$ is the solution set another problem.  Here and  throughout this paper,  we assume that $\Psi\neq\emptyset$. Assume that the set $\Psi$ is a closed  convex subset of a Hilbert space $\mathcal{H}$, then (\ref{e1}) becomes a standard CM problem.
Here we are interested in the case wherein $\Psi$ is the solution
set of another CM problem:
\begin{equation}
\label{e2}
\hbox{minmize}\,\,\{f(x)\,|\, x\in \Omega\}
\end{equation}
i.e., we wish to look at
\begin{equation}
\label{e3}
\Psi:=\{x^*\in \Omega\,|\,f(x^*)\leq f(x)\,|\, \hbox{for all}\,\,x\in \Omega\}
\end{equation}
assuming that $\Psi$ is nonempty. If $f$ is differentiable and we set $F=\nabla f$, then the first order optimality condition of the CM problem (\ref{e2}) translates to the following variational inequality problem of finding a point $x^*\in C$ such that
\begin{equation}
\label{e3}
 \langle F(x^*),x-x^*\rangle\geq0,\quad \forall x\in C.
\end{equation}

The superiorization methodology (SM) strives not to solve (\ref{e1}) but
rather to find a point in $\Psi$ which is superior with respect to $\phi$, i.e., has a lower, but not
necessarily minimal, value of the  objective function $\phi$.
This is done in the SM by first investigating the
bounded perturbation resilience of an algorithm designed to solve (\ref{e2}) and then
proactively using such permitted perturbations in order to steer the iterates of such
an algorithm toward lower values of the $\phi$ objective function while not loosing the
overall convergence to a point in $\Psi$.

So, we aim to prove the bounded perturbation resilience of the PC-algorithms, which will then enable to apply the superiorization idea. To do so, we start by introducing the term \textit{The Basic Algorithm}. Let $\Theta\subseteq \mathcal{H}$ and $\mathfrak{P}$ be any problem with non-empty solution set $\Psi$. Consider the algorithmic operator ${\bf A}_\Psi : \mathcal{H}\rightarrow \Theta$ which works iteratively by
\begin{equation}
\label{BA}
x^{k+1}={\bf A}_\Psi(x^k).
\end{equation}
For any arbitrary starting point $x^0 \in \Theta$. Then (\ref{BA}) is denoted as the \textit{Basic Algorithm}. The bounded perturbation resilience (BPR) of such basic algorithm is defined next.

\begin{definition}
\label{Def41}\cite{HGDC}
 (Bounded perturbation resilience (BPR))
{\rm
An algorithmic operator ${\bf A}_\Psi : \mathcal{H}\rightarrow \Theta$ is said
to be {\it bounded perturbations resilient} if the following is true. If (\ref{BA}) generates sequences $\{x^k\}_{k=0}^\infty$ with $x^0\in \Theta,$ that
converge to points in $\Psi$, then any sequence $\{y^k\}_{k=0}^\infty$, starting from any $y^0\in \Theta,$ generated by
\begin{equation}
\label{mBA}
y^{k+1} = {\bf A_\Psi}(y^k + \lambda_kv^k),\quad \hbox{for all}\,\, k \geq0,
\end{equation}
 also converges to a point in $\Psi$, provided that,
(i) the  sequence $\{v^k\}_{k=0}^\infty$ is bounded, and
(ii) the scalars $\{\lambda_k\}_{k=0}^\infty$
are such that  $\lambda_k\geq 0$ for all
$k \geq 0$, and $\sum_{k=0}^\infty \lambda_k<+\infty$,  and
(iii) $y^k+\lambda_kv^k\in \Theta$ for all $k \geq 0$.
}
\end{definition}

Definition \ref{Def41} is needed only if $\Theta\neq \mathcal{H}$, in which the condition (iii) is enforced in the superiorized version of the basic algorithm, see step (xiv) in the ``Superiorized Version of Algorithm
P" in \cite[p. 5537]{HGDC} and step (14) in ``Superiorized Version of the ML-EM Algorithm"
in \cite[Subsection II.B]{GH}. This will be the case in the present work.

Treating the PC-algorithm as the Basic Algorithm ${\bf A }_ди$ (${\bf A}_\Psi$), our strategy
is to first prove the convergence of Algorithms \ref{PC1-OP} and \ref{PC2-OP} and then show how this yields the BPR of the algorithms according to Definition \ref{Def41}.

A superiorized version of any Basic Algorithm employs the perturbed version of the
Basic Algorithm as in (\ref{mBA}). A certificate to do so in the superiorization method, see
\cite{Censor2}, is gained by showing that the Basic Algorithm is BPR. Therefore, proving the BPR of
an algorithm is the first step toward superiorizing it. This is done for the PC-algorithms
in the next subsection.

\subsection{The BPR of the PC-algorithms}

\qquad
In this subsection, we investigate the bounded perturbation resilience of the PC-algorithms ((\ref{PC-y})-(\ref{PC-x2})).

To this end, we firstly treat the right-hand side of (\ref{PC-x1}) as the algorithmic operator ${\bf A}_\Psi$ of Definition \ref{Def41}, namely, we define for all $k\geq0,$
\begin{equation}
\label{PC1-BP-1}
{\bf A}_\Psi(x^k)=x^k-\gamma \rho_k[(x^k-P_C(x^k-\beta_k F(x^k)))-\beta_k(F(x^k)-F(P_C(x^k-\beta_k F(x^k))))],
\end{equation}
where $\gamma\in(0,2)$,
\begin{equation}
\label{PC1-BP-1a}
\aligned
\beta_k\|F(x^k)-F(P_C(x^k-\beta_k F(x^k)))\|\leq\nu\|x^k-P_C(x^k-\beta_k F(x^k))\|,\,\,\,\, \nu\in(0,1),
\endaligned
\end{equation}
and
\begin{equation}
\label{PC1-BP-2}
\aligned
&\rho_k:=\frac{\| x^k-P_C(x^k-\beta_k F(x^k))\|^2}{\|(x^k-P_C(x^k-\beta_k F(x^k)))-\beta_k (F(x^k)-F(P_C(x^k-\beta_k F(x^k))))\|^2}\\
&\quad-\beta_k\frac{\langle x^k-P_C(x^k-\beta_k F(x^k)), F(x^k)-F(P_C(x^k-\beta_k F(x^k)))\rangle}
{\|(x^k-P_C(x^k-\beta_k F(x^k)))-\beta_k (F(x^k)-F(P_C(x^k-\beta_k F(x^k))))\|^2}
.
\endaligned
\end{equation}

Identify the solution set $\Psi$ with the solution set of the variational inequality problem (\ref{vip}) and identify the additional set $\Theta$ with $C$.

According to Definition \ref{Def41}, we need to show the convergence of any sequence $\{x^k\}_{k=0}^\infty$
 that, starting from any $x^0\in \mathcal{H}$, is generated by
\begin{equation}
\label{PC1-BP-3}
\aligned
x^{k+1}=&x^k+\lambda_kv^k-\gamma \rho_k[(x^k+\lambda_kv^k-P_C(x^k+\lambda_kv^k-\beta_k F(x^k+\lambda_kv^k)))\\
&-\beta_k(F(x^k+\lambda_kv^k)-F(P_C(x^k+\lambda_kv^k-\beta_k F(x^k+\lambda_kv^k))))],\\
\endaligned
\end{equation}
which can be rewritten as follows.

\begin{algorithm} (PC-algorithm I with bounded perturbations)\\
\label{PC1-BP}
Take arbitrarily $x^0\in \mathcal{H}.$
 Given the current iterate $x^k\in \mathcal{H}$, compute
\begin{equation}
\label{PC1-BP-y}
y^k=P_C((x^k+\lambda_kv^k)-\beta_k F((x^k+\lambda_kv^k))),
\end{equation}
where
$\beta_k>0$ is selected to satisfy
\begin{equation}
\label{PC1-BP-beta}
\aligned
\beta_k\|F(x^k+\lambda_kv^k)-F(y^k)\|\leq\nu\|x^k+\lambda_kv^k-y^k\|,\quad \nu\in(0,1).
\endaligned
\end{equation}
Define
\begin{equation}
\label{PC1-BP-d}
d(x^k+\lambda_kv^k,y^k)=(x^k+\lambda_kv^k-y^k)-\beta_k(F(x^k+\lambda_kv^k)-F(y^k)),
\end{equation}
and calculate
\begin{equation}
\label{PC1-BP-x}
x^{k+1}=(x^k+\lambda_kv^k)-\gamma\rho_kd(x^k+\lambda_kv^k,y^k),
\end{equation}
where $\gamma\in(0,2)$,
and
\begin{equation}
\label{PC1-BP-rho}
\aligned
\rho_k:=\frac{\varphi(x^k+\lambda_kv^k,y^k)}
{\|d(x^k+\lambda_kv^k,y^k)\|^2}.
\endaligned
\end{equation}
where $\varphi(x^k+\lambda_kv^k,y^k)=\langle x^k+\lambda_kv^k-y^k,d(x^k+\lambda_kv^k,y^k)\rangle.$
\end{algorithm}

The sequences $\{v^k\}_{k=0}^\infty$ and $\{\lambda_k\}_{k=0}^\infty$ satisfy all the conditions of Definition \ref{Def41}.\vskip 0.2cm

Following the proof of Lemmas \ref{lem31} and \ref{lem32}, we obtain the following lemma.
\begin{lemma}
\label{lem41}
Let $\{\rho_k\}_{k=0}^\infty$ be a sequence which is defined by (\ref{PC1-BP-rho}). Then under Conditions \ref{con:Condition 1.2} and \ref{con:Condition 1.3}, we have
\begin{equation}
\rho_k\geq\frac{1-\nu}{1+\nu^2}.
\end{equation}
\end{lemma}

The next theorem establishes the bounded perturbation resilience of the  PC-algorithm I.
The proof's idea is to build a relationship between BPR and the convergence of
Algorithm \ref{PC1-OP}.

\begin{theorem}
\label{Th41}
Assume that Conditions \ref{con:Condition 1.1}--\ref{con:Condition 1.3} hold. Assume that the sequence $\{v^k\}_{k=0}^\infty$ is bounded, and the positive scalars $\{\lambda_k\}_{k=0}^\infty$ satisfy $\sum_{k=0}^\infty \lambda_k<+\infty$.
Then any sequence $\{x^k\}_{k=0}^\infty$ generated by Algorithm \ref{PC1-BP}  converges weakly to
a solution of the variational inequality problem $(\ref{vip})$.
\end{theorem}

\textit{ Proof.}
Take arbitrarily $x^*\in SOL(C,F)$.
By the definition of $x^{k+1}$, we have
\begin{equation}
\label{th41-1}
\aligned
\|x^{k+1}-x^*\|^2
&=\|x^k+\lambda_kv^k-x^*\|^2+\gamma^2\rho_k^2\|d(x^k+\lambda_kv^k,y^k)\|^2\\
&\quad-2\gamma\rho_k\langle x^k+\lambda_kv^k-x^*,d(x^k+\lambda_kv^k,y^k)\rangle.\\
\endaligned
\end{equation}
Similar with (\ref{th31-2})-(\ref{th31-7}), we have
\begin{equation}
\label{th41-7}
\langle x^k+\lambda_kv^k-x^*,d(x^k+\lambda_kv^k,y^k)\rangle\geq\varphi(x^k+\lambda_kv^k,y^k).
\end{equation}
Substituting (\ref{th41-7}) into (\ref{th41-1}) and using $\rho_k=\varphi(x^k+\lambda_kv^k,y^k)/\|d(x^k+\lambda_kv^k,y^k)\|^2$, we have
\begin{equation}
\label{th41-8}
\aligned
&\|x^{k+1}-x^*\|^2\\
&\leq\|x^k+\lambda_kv^k-x^*\|^2+\gamma^2
\rho_k^2\|d(x^k+\lambda_kv^k,y^k)\|^2\\
&\quad-2\gamma\rho_k\varphi(x^k+\lambda_kv^k,y^k)\\
&=\|x^k+\lambda_kv^k-x^*\|^2-\gamma(2-\gamma)\rho_k\varphi(x^k+\lambda_kv^k,y^k).\\
\endaligned
\end{equation}
Again, using the definition of $x^{k+1}$, we have
\begin{equation}
\label{th41-9}
\aligned
\rho_k\varphi(x^k+\lambda_kv^k,y^k)&=\|\rho_kd(x^k+\lambda_kv^k,y^k)\|^2\\
&=\frac{1}{\gamma^2}\|x^{k+1}-(x^k+\lambda_kv^k)\|^2.\\
\endaligned
\end{equation}
Combining the inequalities (\ref{th41-8}) and (\ref{th41-9}), we obtain
\begin{equation}
\label{th41-10}
\aligned
\|x^{k+1}-x^*\|^2\leq\|x^k+\lambda_kv^k-x^*\|^2-\frac{2-\gamma}{\gamma}\|x^{k+1}-(x^k+\lambda_kv^k)\|^2.
\endaligned
\end{equation}
By Lemma \ref{lem-ineq}, we have
\begin{equation}
\label{th41-11}
\aligned
\|x^k+\lambda_kv^k-x^*\|^2&=\|x^k-x^*\|^2+\lambda_k^2\|v^k\|^2+2\lambda_k\langle x^k-x^*,v^k\rangle\\
&\leq(1+\lambda_k\|v^k\|)\|x^k-x^*\|^2+\lambda_k^2\|v^k\|^2+\lambda_k\|v^k\|,\\
\endaligned
\end{equation}
and
\begin{equation}
\label{th41-12}
\aligned
&\|x^{k+1}-(x^k+\lambda_kv^k)\|^2\\
&=\|x^{k+1}-x^k\|^2+\lambda_k^2\|v^k\|^2-2\lambda_k\langle x^{k+1}-x^k,v^k\rangle\\
&\geq(1-\lambda_k\|v^k\|)\|x^{k+1}-x^k\|^2+\lambda_k^2\|v^k\|^2-\lambda_k\|v^k\|.\\
\endaligned
\end{equation}
Combining (\ref{th41-10})-(\ref{th41-12}), we obtain
\begin{equation}
\label{th41-13}
\aligned
&\|x^{k+1}-x^*\|^2\\
&\leq(1+\lambda_k\|v^k\|)\|x^k-x^*\|^2+
\frac{2\lambda_k}\gamma\|v^k\|-\frac{2(1-\gamma)}\gamma\lambda_k^2\|v^k\|^2\\
&\quad-\frac{2-\gamma}{\gamma}(1-\lambda_k\|v^k\|)\|x^{k+1}-x^k\|^2.\\
\endaligned
\end{equation}
From the assumptions on $\{\lambda_k\}_{k=0}^\infty$ and the fact that $\{v^k\}_{k=0}^\infty$ is bounded, we have
\begin{equation}
\label{th41-14}
\sum_{k=0}^\infty \lambda_k\|v^k\|<+\infty,\quad \sum_{k=0}^\infty \lambda_k^2\|v^k\|^2<+\infty,
\end{equation}
which means that
\begin{equation}
\label{th41-15}
\lim_{k\rightarrow\infty}\lambda_k\|v^k\|=0,\quad \lim_{k\rightarrow\infty}\lambda_k^2\|v^k\|^2=0.
\end{equation}
Assume that $\lambda_k\|v^k\|\in[0,\mu)$, where $\mu\in[0,1)$, then we get
\begin{equation}
\label{th41-16}
\aligned
&\|x^{k+1}-x^*\|^2\\
&\leq(1+\lambda_k\|v^k\|)\|x^k-x^*\|^2+
\frac{2\lambda_k}\gamma\|v^k\|-\frac{2(1-\gamma)}\gamma\lambda_k^2\|v^k\|^2\\
&\quad-\frac{2-\gamma}{\gamma}(1-\mu)\|x^{k+1}-x^k\|^2\\
&\leq(1+\lambda_k\|v^k\|)\|x^k-x^*\|^2+
\frac{2\lambda_k}\gamma\|v^k\|-\frac{2(1-\gamma)}\gamma\lambda_k^2\|v^k\|^2.
\endaligned
\end{equation}
Following the proof of (\ref{th31-15}), we get
\begin{equation}
\label{th41-17}
\aligned
\lim_{k\rightarrow\infty}\|x^{k+1}-x^k\|=0.
\endaligned
\end{equation}
Similar to (\ref{th31-16}), we get
\begin{equation}
\label{th41-18}
\aligned
\varphi(x^k+\lambda_kv^k,y^k)\leq\frac{2(1+\nu^2)}{(1-\nu)\gamma^2}[\|x^{k+1}-x^k\|^2+\lambda_k^2\|v^k\|^2].
\endaligned
\end{equation}
Using Lemma \ref{lem-ineq} and the proof of (\ref{le12}), we have
\begin{equation}
\label{th41-19}
\aligned
\varphi(x^k+&\lambda_kv^k,y^k)\\
&\geq(1-\nu)\|x^k+\lambda_kv^k-y^k\|^2\\
&\geq(1-\nu)[(1-\lambda_k\|v^k\|)\|x^k-y^k\|^2+\lambda_k^2\|v^k\|^2-\lambda_k\|v^k\|]\\
&\geq(1-\nu)[(1-\mu)\|x^k-y^k\|^2+\lambda_k^2\|v^k\|^2-\lambda_k\|v^k\|].\\
\endaligned
\end{equation}
From (\ref{th41-18}) and (\ref{th41-19}), we obtain
\begin{equation}
\label{th41-20}
\aligned
\|x^k-y^k\|^2&\leq\frac1{(1-\mu)}
\left(\frac{2(1+\nu^2)}{[(1-\nu)\gamma]^2}[\|x^{k+1}-x^k\|^2+\lambda_k^2\|v^k\|^2]+\lambda_k\|v^k\|\right),\\
\endaligned
\end{equation}
which with (\ref{th41-15}) and (\ref{th41-17}) yields
\begin{equation}
\label{th41-21}
\aligned
\lim_{k\rightarrow\infty}\|x^k-y^k\|^2=0.
\endaligned
\end{equation}
Now following the lines of Theorem \ref{th31} the rest of the proof is completed. $\Box$
\vskip 0.2cm

Similar to Theorem \ref{Th32}, we get the convergence rate of Algorithm \ref{PC1-BP}.
\begin{theorem}
\label{Th46} Assume that Conditions \ref{con:Condition 1.1}--\ref{con:Condition 1.3} hold. Let $%
\{x^{k}\}_{k=0}^{\infty }$ and $\{y^{k}\}_{k=0}^{\infty }$ be any two sequences generated by
Algorithm \ref{PC1-BP}. For any integer $t>0$, we have a $y_{t}\in C$
which satisfies
\begin{equation}
\label{Th46-0a}
\aligned
\langle F(x), y_{t}-x\rangle \leq \frac{1}{2\gamma \Upsilon _{t}}%
(\Vert x-x^{0}\Vert ^{2}+2M),\quad \forall x\in C,\endaligned
\end{equation}%
where
\begin{equation}
\label{Th46-0b}
 y_{t}=\frac{1}{\Upsilon _{t}}\sum_{k=0}^{t}\rho_k\beta_ky^{k},\,\, \Upsilon
_{t}=\sum_{k=0}^{t}\rho_k\beta_k,\,\, \hbox{and}\,\, M=\sup_{k\in {\Bbb N}}\{\Vert x^{k}-x\Vert
\}\sum_{k=0}^{\infty }\lambda_k\Vert v^k\Vert.
\end{equation}
\end{theorem}

Next, we investigate the bounded perturbation resilience of the PC-algorithm II.
We treat the right-hand side of (\ref{PC-x2}) as the algorithmic operator ${\bf A}_\Psi$ of Definition \ref{Def41}, namely, we define for all $k\geq0,$
\begin{equation}
{\bf A}_\Psi(x^k)=P_C[x^k-\gamma \beta_k\rho_kF(P_C(x^k-\beta_k F(x^k)))],
\end{equation}
where $\gamma\in(0,2)$, $\beta_k$ and $\rho_k$ are defined as in (\ref{PC1-BP-1a}) and (\ref{PC1-BP-2}), respectively.

According to Definition \ref{Def41}, we need to show the convergence of any sequence $\{x^k\}_{k=0}^\infty$
generated by
\begin{equation}
\aligned
x^{k+1}=P_C[x^k+\lambda_kv^k-\gamma \beta_k\rho_kF(P_C(x^k+\lambda_kv^k-\beta_k F(x^k+\lambda_kv^k)))],
\endaligned
\end{equation}
for any starting point $x^0\in \mathcal{H}$.

\begin{algorithm} (PC-algorithm II with bounded perturbations)\\
\label{PC2-BP}
Take arbitrarily $x^0\in \mathcal{H}.$
 Given the current iterate $x^k\in \mathcal{H}$, compute
\begin{equation}
\label{PC1-BP-y}
y^k=P_C((x^k+\lambda_kv^k)-\beta_k F((x^k+\lambda_kv^k))),
\end{equation}
where $\beta_k>0$ is selected via (\ref{PC1-BP-beta}).

Define $d(x^k+\lambda_kv^k,y^k)$ as in (\ref{PC1-BP-d}).
Calculate
\begin{equation}
\label{PC2-BP-x}
x^{k+1}=P_C[x^k+\lambda_kv^k-\gamma \beta_k\rho_kF(y^k)],
\end{equation}
where $\gamma\in(0,2)$, $\rho_k$ is defined as in (\ref{PC1-BP-rho}).
\end{algorithm}

The sequence $\{v^k\}_{k=0}^\infty$ and the scalars $\{\lambda_k\}_{k=0}^\infty$ satisfy all the conditions in Definition \ref{Def41}.\vskip 0.1cm

Following the proof of Theorems \ref{th33} and \ref{Th41}, we get the convergence of Algorithm \ref{PC2-BP}.
\begin{theorem}
\label{Th42}
Assume that Conditions \ref{con:Condition 1.1}--\ref{con:Condition 1.3} hold. Assume that the sequence $\{v^k\}_{k=0}^\infty$ is bounded, and the positive scalars $\{\lambda_k\}_{k=0}^\infty$ satisfy $\sum_{k=0}^\infty \lambda_k<+\infty$.
Then any sequence $\{x^k\}_{k=0}^\infty$   generated by Algorithm \ref{PC2-BP}  converges weakly to
a solution of the variational inequality problem $(\ref{vip})$.
\end{theorem}

\begin{theorem}
\label{Th49} Assume that Conditions \ref{con:Condition 1.1}--\ref{con:Condition 1.3} hold. Let $%
\{x^{k}\}_{k=0}^{\infty }$ and $\{y^{k}\}_{k=0}^{\infty }$ be any two sequences generated by
Algorithm \ref{PC1-BP}. For any integer $t>0$, we have a $y_{t}\in C$
which satisfies
\begin{equation}
\label{Th49-1}
\aligned
\langle F(x), y_{t}-x\rangle \leq \frac{1}{2\gamma \Upsilon _{t}}%
(\Vert x-x^{0}\Vert ^{2}+2M),\quad \forall x\in C,\endaligned
\end{equation}%
where
\begin{equation}
\label{Th49-2}
 y_{t}=\frac{1}{\Upsilon _{t}}\sum_{k=0}^{t}\rho_k\beta_ky^{k},\,\, \Upsilon
_{t}=\sum_{k=0}^{t}\rho_k\beta_k,\,\, \hbox{and}\,\, M=\sup_{k\in {\Bbb N}}\{\Vert x^{k+1}-x\Vert
\}\sum_{k=0}^{\infty }\lambda_k\Vert v^{k}\Vert.
\end{equation}

\end{theorem}

\section{Construction of the inertial PC-algorithms}\label{sec:iner-PC}

\qquad
In this section, we construct four classes of inertial PC-algorithms by using outer perturbations and bounded perturbations, i.e., identifying the $e_i^k$, $k=1,2$ and $\lambda_k$, $v^k$ with special values.

The inertial-type algorithms originate from the \textit{heavy ball method}
of the second-order dynamical systems in time \cite{Alvarez} and speed up the original algorithm without the inertial effects.
Recently there are increasing interests in studying inertial-type algorithms, see for example \cite{Alvarez,APR,OBP,Dong3} and the references therein.

Using Algorithm \ref{PC1-OP}, we construct the following inertial PC-algorithm I (iPC I-1 for short):
\begin{equation}
\label{IPC1-OP}
\left\{
\aligned
&y^k=P_C(x^k-\beta_k F(x^k))+\alpha_k^{(1)}(x^k-x^{k-1}),\\
&d(x^k,y^k-\alpha_k^{(1)}(x^k-x^{k-1}))=(x^k-y^k+\alpha_k^{(1)}(x^k-x^{k-1}))\\
&\quad-\beta_k(F(x^k)-F(y^k-\alpha_k^{(1)}(x^k-x^{k-1}))),\\
&x^{k+1}=x^k-\gamma \rho_kd(x^k,y^k-\alpha_k^{(1)}(x^k-x^{k-1}))+\alpha_k^{(2)}(x^k-x^{k-1}).
\endaligned
\right.
\end{equation}
where
$\gamma\in(0,2)$,  $\alpha _{k}^{(i)}\in[0,1]$, $i=1,2$, and
 $\beta_k>0$ is selected to satisfy
\begin{equation}
\label{IPC1-OP-beta}
\aligned
\beta_k\|F(x^k)-&F(y^k-\alpha_k^{(1)}(x^k-x^{k-1})\|\\
&\leq\nu\|x^k-y^k+\alpha_k^{(1)}(x^k-x^{k-1}\|,\,\,\,\, \nu\in(0,1),
\endaligned
\end{equation}
and
\begin{equation}
\label{IPC1-OP-rho}
\rho_k:=\frac{\langle x^k-y^k+\alpha_k^{(1)}(x^k-x^{k-1})),d(x^k,y^k-\alpha_k^{(1)}(x^k-x^{k-1}))\rangle}
{\|d(x^k,y^k-\alpha_k^{(1)}(x^k-x^{k-1}))\|^2}.
\end{equation}

For the convergence of the inertial algorithm, the following condition should be imposed on the inertial
parameters $\alpha _{k}^{(i)}$, $i=1,2,$
\begin{equation}
\label{IPC1-OP-alhpa}
\sum_{k=0}^\infty\alpha _{k}^{(i)}\|x^k-x^{k-1}\|<+\infty,\quad i=1,2.
\end{equation}

\begin{remark}
\textrm{
Condition (\ref{IPC1-OP-alhpa}) can be enforced by a simple online updating rule such as, given $\alpha^{(i)}\in [0,1]$, $i=1,2$
\begin{equation}
\label{N2}
\alpha_k^{(i)}=\min\{\alpha^{(i)},\zeta_k^{(i)}\},
\end{equation}
 where $\zeta_k^{(i)}>0$,
$\zeta_k^{(i)}\|x^k- x^{k-1}\|$ is summable. For instance, one can choose
\begin{equation}
\zeta_k^{(i)}=\frac{\zeta^{(i)}}{k^{1+\xi}\|x^k- x^{k-1}\|},\quad \zeta^{(i)}>0,\quad\xi>0.
\end{equation}
In practical calculation, $\|x^k- x^{k-1}\|$ rapidly vanishes as $k\rightarrow\infty$. So most of the time,
with proper choice of $\alpha^{(i)}$, (\ref{N2}) may never be triggered.
}
\end{remark}
\vskip 0.2cm

Similar with Theorem \ref{th31}, we get the convergence of the inertial PC-algorithm I (\ref{IPC1-OP}).
\begin{theorem}
\label{Th51}
Assume that Conditions \ref{con:Condition 1.1}--\ref{con:Condition 1.3} hold. Assume that the sequences $\{\alpha_k^{(i)}\}_{k=0}^\infty$, $i=1,2$ satisfy  (\ref{IPC1-OP-alhpa}). Then any sequence $\{x^k\}_{k=0}^\infty$ generated by the inertial PC-algorithm I (\ref{IPC1-OP})   converges weakly to
a solution of the variational inequality problem $(\ref{vip})$.
\end{theorem}

Using Algorithm \ref{PC2-OP}, we construct the following inertial PC-algorithm II (iPC II-1):
\begin{equation}
\label{IPC2-OP}
\left\{
\aligned
&y^k=P_C(x^k-\beta_k F(x^k)+\alpha_k^{(1)}(x^k-x^{k-1})),\\
&x^{k+1}=P_C(x^k-\gamma \beta_k\rho_k F(y^k)+\alpha_k^{(2)}(x^k-x^{k-1})).
\endaligned
\right.
\end{equation}
where $\gamma\in(0,2)$,  $\alpha _{k}^{(i)}\in[0,1]$, $i=1,2$ and $\beta_k>0$ is selected to satisfy
\begin{equation}
\label{IPC2-OP-beta}
\aligned
\beta_k\|F(x^k)-F(y^k)\|\leq\nu\|x^k-y^k\|,\,\,\,\, \nu\in(0,1),
\endaligned
\end{equation}
and
\begin{equation}
\label{IPC2-OP-rho}
\aligned
\rho_k:=\frac{\langle x^k-y^k, d(x^k,y^k)\rangle}{\|d(x^k,y^k)\|^2},
\endaligned
\end{equation}
and
\begin{equation}
\label{IPC2-OP-d}
d(x^k,y^k):=(x^k-y^k)-\beta_k(F(x^k)-F(y^k))+\alpha_k^{(1)}(x^k-x^{k-1}).
\end{equation}

\vskip 0.2cm

\begin{theorem}
\label{Th52}
Assume that Conditions \ref{con:Condition 1.1}--\ref{con:Condition 1.3} hold. Assume that the sequences $\{\alpha_k^{(i)}\}_{k=0}^\infty$, $i=1,2$ satisfy  (\ref{IPC1-OP-alhpa}) and
 \begin{equation}
\label{Th52-1}
\alpha_k^{(1)}\|x^{k}-x^{k-1}\|\leq\mu\|x^k-y^k\|,
\end{equation}
where $\mu\in [0,\nu).$
Then any sequence $\{x^k\}_{k=0}^\infty$ generated by the inertial PC-algorithm II (\ref{IPC2-OP})  converges weakly to
a solution of the variational inequality problem $(\ref{vip})$.
\end{theorem}

Using Algorithm \ref{PC1-BP},
we construct the following inertial PC-algorithm I (iPC I-2):
 \begin{equation}
\label{IPC1-BP}
\left\{
\aligned
&w^k=x^k+\alpha_k(x^k-x^{k-1})\\
&y^k=P_C(w^k-\beta_k F(w^k)),\\
&d(w^k,y^k)=(w^k-y^k)-\beta_k(F(w^k)-F(y^k)),\\
&x^{k+1}=w^k-\gamma\rho_kd(w^k,y^k).
\endaligned
\right.
\end{equation}
where
$\gamma\in(0,2)$,  $\alpha _{k}^{(i)}\in[0,1]$, $i=1,2$ and $\beta_k>0$ is selected to satisfy
\begin{equation}
\label{IPC1-BP-beta}
\aligned
\beta_k\|F(w^k)-F(y^k)\|\leq\nu\|w^k-y^k\|,\,\,\,\, \nu\in(0,1),
\endaligned
\end{equation}
and
\begin{equation}
\label{IPC1-BP-rho}
\rho_k:=\frac{\langle w^k-y^k,d(w^k,y^k)\rangle}
{\|d(w^k,y^k)\|^2}.
\end{equation}
\vskip 0.2cm

We extend Theorem \ref{Th41} to the convergence of the inertial PC-algorithm II.
\begin{theorem}
\label{Th53}
Assume that Conditions \ref{con:Condition 1.1}--\ref{con:Condition 1.3} hold. Assume that the sequence $\{\alpha_k^{(i)}\}_{k=0}^\infty$, $i=1,2$ satisfy  (\ref{IPC1-OP-alhpa}). Then any sequence $\{x^k\}_{k=0}^\infty$ generated by the inertial PC-algorithm I (\ref{IPC1-BP})  converges weakly to
a solution of the variational inequality problem $(\ref{vip})$.

\end{theorem}

Using Algorithm \ref{PC2-BP},
we construct the following inertial PC-algorithm II (iPC II-2):
 \begin{equation}
\label{IPC2-BP}
\left\{
\aligned
&w^k=x^k+\alpha_k(x^k-x^{k-1})\\
&y^k=P_C(w^k-\beta_k F(w^k)),\\
&x^{k+1}=P_C(w^k-\gamma\beta_k\rho_k F(y^k)).
\endaligned
\right.
\end{equation}
where $\gamma\in(0,2)$, $\alpha _{k}^{(i)}\in[0,1]$, $i=1,2$, and $\beta_k$ and $\rho_k$ are defined as
 (\ref{IPC1-BP-beta}) and (\ref{IPC1-BP-rho}), respectively.
\vskip 0.2cm

We extend Theorem \ref{Th49} to the convergence of the inertial PC-algorithms II.
\begin{theorem}
\label{Th54}
Assume that Conditions \ref{con:Condition 1.1}--\ref{con:Condition 1.3} hold. Assume that the sequence $\{\alpha_k^{(i)}\}_{k=0}^\infty$, $i=1,2$ satisfy  (\ref{IPC1-OP-alhpa}). Then any sequence $\{x^k\}_{k=0}^\infty$ generated by the inertial PC-algorithm II (\ref{IPC2-BP})  converges weakly to
a solution of the variational inequality problem $(\ref{vip})$.

\end{theorem}

\begin{remark}
\label{Rem56}
In \cite{Dong2}, by using a different technique, the authors proved the convergence of the inertial PC-algorithm I (\ref{IPC1-BP}) provided that $\{\alpha_k\}_{k=0}^\infty$ is nondecreasing with $\alpha_1=0,$ $0\leq\alpha_k\leq\alpha<1,$ and $\sigma,\delta>0$ are such that
 \begin{equation}
\label{Rem56a}
\delta>\frac{\alpha^2(1+\alpha)+\alpha\sigma}{1-\alpha^2},\quad
0<\gamma\leq\frac{2\left[\delta-\alpha((1+\alpha)+\alpha\delta+\sigma)\right]}
{\delta[1+\alpha(1+\alpha)+\alpha\delta+\sigma]}.
\end{equation}
They showed  the efficiency and advantage of the inertial PC-algorithm I (\ref{IPC1-BP}) with above inertial parameters through numerical experiments. But, inertial variants of the PC-algorithm II was not considered in \cite{Dong2}!
\end{remark}

\section{Numerical experiments\label{sec:num}}

\qquad In this section, we compare and illustrate the performances of all the presented algorithms for the problem of sparse signal recovery problem. The algorithms are: the PC-algorithm I (\ref{PC-x1}) (PC I), the PC-algorithm II (\ref{PC-x2}) (PC II) the inertial PC-algorithm I (\ref{IPC1-OP}) (iPC I-1),
the inertial PC-algorithm II (\ref{IPC2-OP}) (iPC II-1), the inertial PC-algorithm I (\ref{IPC1-BP}) (iPC I-2), the inertial PC-algorithm I (\ref{IPC2-BP}) (iPC II-2)
and the inertial PC-algorithm I (\ref{IPC1-BP}) with the inertial parameters satisfying the conditions in Remark \ref{Rem56} (iPC I for short). \vskip 2mm

Choose the following set of parameters. Take $\sigma =5$, $\rho =0.9$ , $\mu =0.7$ and $\gamma=1$. For iPC I-1 and iPC II-1,  set

\begin{equation}
\alpha_k^{(i)}=\min\{\alpha^{(i)},\zeta_k^{(i)}\},
\end{equation}
 where $\alpha^{(i)}\in(0,1)$, and
\begin{equation}
\zeta_k^{(i)}=\frac{\alpha^{(i)}}{k^2\|x^k- x^{k-1}\|}.
\end{equation}
Similarly,
for iPC I-2 and iPC II-2,  set
\begin{equation}
\alpha_k=\min\{\alpha,\zeta_k\},
\end{equation}
 where $\alpha\in(0,1)$, and
\begin{equation}
\zeta_k=\frac{\alpha}{k^2\|x^k- x^{k-1}\|}.
\end{equation}
Take $\alpha^{(i)}=0.4$
in iPC I-1.
In order to guarantee the convergence of iPC II-1, the inertial parameters $\alpha_k^{(1)}$  should satisfy the condition (\ref{Th52-1}). After running numerous simulations, we find that condition (\ref{Th52-1}) is satisfied when $\alpha^{(1)}$ is taken in $(0,0.4]$. So, we decided to choose $\alpha^{(i)}=0.4$ in the presented example. We also take $\alpha_k=0.8$ for iPC I-2 and iPC II-2, and  $\alpha_k=0.79$ for  iPC I, respectively.

\begin{figure}[htp]
\setlength{\floatsep}{0pt} \setlength{\abovecaptionskip}{-10pt}
\centering
\scalebox{0.6} {\includegraphics{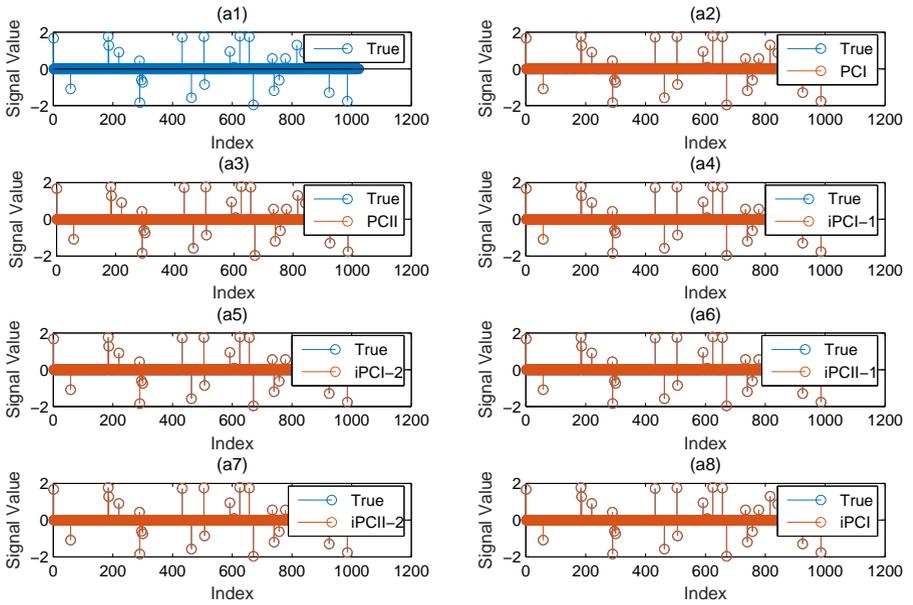}}
\caption{(a1) is the true sparse signal, (a2)-(a8) are the recovered signal vs the true signal by "PC I", "PC II", "iPC I-1", "iPC I-2", "iPC II-1" "iPC II-2" and "iPC I", respectively. }\label{ex1_k30}
\end{figure}

\begin{figure}[htp]
\setlength{\floatsep}{0pt} \setlength{\abovecaptionskip}{-5pt}
\centering
\scalebox{0.5} {\includegraphics{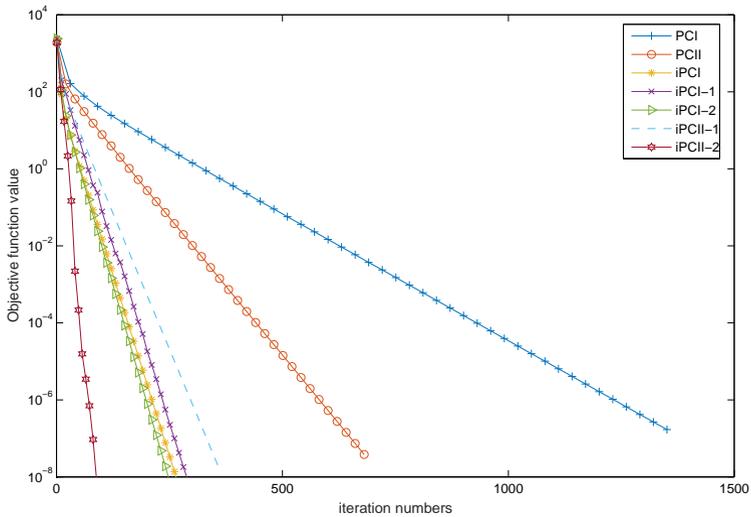}}
\caption{Comparison of the objective function value versus the iteration numbers of the different methods. }\label{ex1_k30_fun}
\end{figure}

\begin{example}
Let $x^0 \in \mathbb{R}^{n}$ be a $K$-sparse signal, $K\ll n$. The sampling matrix $A \in \mathbb{R}^{m\times n} (m<<n)$ is stimulated by standard Gaussian distribution and vector $b = Ax^0 + e$, where $e$ is additive noise. When $e=0$, it means that there is no noise to the observed data. Our task is to recover the signal $x^0$ from the data $b$.
\end{example}
It's well-known that the sparse signal $x^0$ can be recovered by solving the following LASSO problem \cite{Tibshirani96},
\begin{equation}
\begin{aligned}
\min_{x\in \mathbb{R}^n}\ & \frac{1}{2}\| Ax - b \|_{2}^{2} \\
s.t. \ & \|x\|_1 \leq t,
\end{aligned}\label{ex:sparse}
\end{equation}
where $t >0$. It is easy to see that the optimization problem (\ref{ex:sparse}) is a special case of the variational inequality problem (\ref{vip}), where $F(x) = A^{T}(Ax-b)$ and $C = \{ x \mid \|x\|_1 \leq t \}$. We can use the proposed iterative algorithms to solve the optimization problem (\ref{ex:sparse}). Although the orthogonal projection onto the closed convex set $C$ doesn't have a closed-form solution, the projection operator $P_{C}$ can be precisely computed in a polynomial time (see for example \cite{dssc08}).
The following inequality was defined as the stopping criteria,
\begin{equation}
\| x^{k+1} - x^k \| \leq \epsilon,
\end{equation}
where $\epsilon >0$ is a given small constant. $``Iter"$ denotes the iteration numbers. $``Obj"$ represents the objective function value and $``Err"$ is the $2$-norm error between the recovered signal and the true $K$-sparse signal. We divide the experiments into two parts. One task is to recover the sparse signal $x^0$ from noise observation vector $b$ and the other is to recover the sparse signal from noiseless data $b$. For the noiseless case, the obtained numerical results are reported in Table \ref{example1:results1}. To visually view the results, Figure \ref{ex1_k30} shows the recovered signal compared with the true signal $x_0$ when $K=30$. We can see from Figure \ref{ex1_k30} that the recovered signal is the same as the true signal. Further, Figure \ref{ex1_k30_fun} presents the objective function value versus the iteration numbers.

\begin{table}[htp]
\caption{Numerical results obtained by the proposed iterative algorithms
when $m =240, n =1024$ in the noiseless case.}
\centering
\resizebox{\columnwidth}{!}{
\begin{tabular}{c|c|cccccccc}
\hline
\multirow{2}[1]{*}{$K$-sparse\vspace{.1in}}  & \multirow{2}[1]{*}{Methods} & &   \multicolumn{3}{c}{$ \epsilon = 10^{-4}$} &  & \multicolumn{3}{c}{$ \epsilon = 10^{-6}$} \\
\cline{4-6} \cline{8-10}
signal & & & $Iter$ & $Obj$ & $Err$ & &  $Iter$ &  $Obj$ & $Err$  \\
\hline
 \hline
\multirow{7}[1]{*}{$K=20$}
& PC I & & $443$ & $8.5443e-4$ & $0.0076$ & & $841$ & $8.8686e-8$ & $8.0044e-5$\\
& PC II & & $227$ & $1.4410e-4$ & $0.0032$ & & $409$ & $1.6627e-8$ & $3.5157e-5$ \\
& iPC I & & $72$ & $7.4506e-6$ & $7.0980e-4$ & & $103$ & $1.0550e-9$ & $8.6892e-6$ \\
& iPC I-1 & & $118$ & $2.4372e-5$ & $0.0013$ & & $191$ & $1.6450e-9$ & $1.1226e-5$ \\
& iPC I-2 & & $56$ & $6.0246e-6$ & $3.8709e-4$ & & $99$ & $7.6703e-10$ & $7.3893e-6$ \\
& iPC II-1 & & $131$ & $3.5673e-5$ & $0.0016$ & & $231$ & $4.4078e-9$ & $1.7965e-5$ \\
& iPC II-2 & & $55$ & $1.2499e-6$ & $1.3129e-4$ & & $80$ & $1.6088e-10$ & $3.0994e-6$ \\
\hline
\multirow{7}[1]{*}{$K=30$}
& PC I & & $732$ & $0.0019$ & $0.0166$ & & $1361$ & $1.4492e-7$ & $1.3025e-4$\\
& PC II & & $407$ & $3.8218e-4$ & $0.0076$ & & $689$ & $2.9161e-8$ & $5.9125e-5$ \\
& iPC I & & $178$ & $5.6970e-5$ & $0.0029$ & & $275$ & $4.0179e-9$ & $2.1673e-5$ \\
& iPC I-1 & & $204$ & $4.3636e-5$ & $0.0026$ & & $295$ & $5.7110e-9$ & $2.5667e-5$ \\
& iPC I-2 & & $170$ & $5.2038e-5$ & $0.0028$ & & $259$ & $3.4741e-9$ & $2.0260e-5$ \\
& iPC II-1 & & $235$ & $9.0708e-5$ & $0.0037$ & & $373$ & $7.1864e-9$ & $2.9215e-5$ \\
& iPC II-2 & & $73$ & $3.9273e-6$ & $7.7445e-4$ & & $100$ & $4.9113e-11$ & $2.1110e-6$ \\
\hline
\multirow{7}[1]{*}{$K=40$}
& PC I & & $1440$ & $0.0041$ & $0.0376$ & & $4548$ & $5.7806e-7$ & $5.2196e-4$\\
& PC II & & $867$ & $0.0010$ & $0.0188$ & & $2459$ & $1.4600e-7$ & $2.6263e-4$ \\
& iPC I & & $398$ & $1.6322e-4$ & $0.0075$ & & $1074$ & $2.2549e-8$ & $1.0309e-4$ \\
& iPC I-1 & & $362$ & $4.6840e-4$ & $0.0116$ & & $1006$ & $2.3931e-8$ & $1.0636e-4$ \\
& iPC I-2 & & $382$ & $1.4364e-4$ & $0.0071$ & & $1009$ & $1.9430e-8$ & $9.5612e-5$ \\
& iPC II-1 & & $498$ & $2.5750e-4$ & $0.0095$ & & $1329$ & $3.5183e-8$ & $1.2876e-4$ \\
& iPC II-2 & & $210$ & $3.1377e-5$ & $0.0033$ & & $520$ & $4.3266e-9$ & $4.5213e-5$ \\
 \hline
\hline
\end{tabular}\label{example1:results1}
}
\end{table}

For the noise observation $b$, we assume that the vector $e$ is corrupted by Gaussian noise with zero mean and $\beta$ variances. The system matrix $A$ is the same as the noiseless case and the sparsity level $K=30$.
We list the numerical results for different noise level $\beta$ in Table \ref{example1:results2}.

\begin{table}[htp]
\centering
\caption{Numerical results for the proposed iterative algorithms with different noise value $\beta$.}
\begin{tabular}{c|c|cccccccc}
\hline
\multirow{2}[1]{*}{Variances\vspace{1.in}}  & \multirow{2}[1]{*}{Methods} & &   \multicolumn{3}{c}{$ \epsilon = 10^{-4}$} &  & \multicolumn{3}{c}{$ \epsilon = 10^{-6}$} \\ \cline{4-6} \cline{8-10}
$\beta = $ & & & $Iter$ & $Obj$ & $Err$ & &  $Iter$ &  $Obj$ & $Err$  \\
\hline
 \hline
\multirow{7}[1]{*}{$0.01$}
& PC I & & $602$ & $0.0065$ & $0.0158$ & & $1435$ & $0.0050$ & $0.0130$\\
& PC II & & $314$ & $0.0058$ & $0.0124$ & & $739$ & $0.0050$ & $0.0130$ \\
& iPC I & & $118$ & $0.0057$ & $0.0111$ & & $317$ & $0.0050$ & $0.0130$ \\
& iPC I-1 & & $146$ & $0.0057$ & $0.0114$ & & $342$ & $0.0050$ & $0.0130$ \\
& iPC I-2 & & $110$ & $0.0057$ & $0.0110$ & & $301$ & $0.0050$ & $0.0130$ \\
& iPC II-1 & & $170$ & $0.0057$ & $0.0115$ & & $402$ & $0.0050$ & $0.0130$ \\
& iPC II-2 & & $58$ & $0.0057$ & $0.0109$ & & $137$ & $0.0050$ & $0.0130$ \\
\hline
\multirow{7}[1]{*}{$0.02$}
& PC I & & $743$ & $0.0213$ & $0.0308$ & & $954$ & $0.0224$ & $0.0185$\\
& PC II & & $395$ & $0.0203$ & $0.0262$ & & $466$ & $0.0224$ & $0.0185$ \\
& iPC I & & $161$ & $0.0200$ & $0.0242$ & & $175$ & $0.0224$ & $0.0185$ \\
& iPC I-1 & & $181$ & $0.0201$ & $0.0245$ & & $222$ & $0.0224$ & $0.0185$ \\
& iPC I-2 & & $152$ & $0.0200$ & $0.0241$ & & $157$ & $0.0224$ & $0.0185$ \\
& iPC II-1 & & $216$ & $0.0201$ & $0.0247$ & & $254$ & $0.0224$ & $0.0185$ \\
& iPC II-2 & & $65$ & $0.0200$ & $0.0235$ & & $89$ & $0.0224$ & $0.0185$ \\
\hline
\multirow{7}[1]{*}{$0.05$}
& PC I & & $564$ & $0.1439$ & $0.0541$ & & $1400$ & $0.1344$ & $0.0608$\\
& PC II & & $302$ & $0.1432$ & $0.0522$ & & $683$ & $0.1344$ & $0.0607$ \\
& iPC I & & $112$ & $0.1430$ & $0.0511$ & & $298$ & $0.1344$ & $0.0607$ \\
& iPC I-1 & & $181$ & $0.1430$ & $0.0515$ & & $318$ & $0.1344$ & $0.0607$ \\
& iPC I-2 & & $105$ & $0.1430$ & $0.0511$ & & $281$ & $0.1344$ & $0.0607$ \\
& iPC II-1 & & $171$ & $0.1430$ & $0.0515$ & & $371$ & $0.1344$ & $0.0607$ \\
& iPC II-2 & & $66$ & $0.1430$ & $0.0510$ & & $110$ & $0.1344$ & $0.0607$ \\

 \hline
\hline

\end{tabular}\label{example1:results2}
\end{table}

From Tables \ref{example1:results1} and \ref{example1:results2}, and Figure \ref{ex1_k30}, we conclude: (i) PC II behaves better than PC I;
(ii) the inertial type algorithms improve the original algorithms; (iii)  iPC II-2
has best performance among the inertial type algorithms, while iPC II-1 behaves worst; (iii)
the performance of iPC1, iPC1-1 and iPC1-2 is close and almost same.

\newpage
\section*{Acknowledgement}
We wish to thank the anonymous referees for the thorough analysis and review, their comments and suggestions helped tremendously in improving the quality of this paper and made it suitable for publication.

The first author is supported by National Natural Science Foundation of China (No. 71602144) and
Open Fund of Tianjin Key Lab for Advanced Signal Processing (No. 2016ASP-TJ01). The second author is supported by the EU FP7 IRSES program STREVCOMS (No. PIRSES-GA-2013-612669).


\end{document}